\ifx\documentclass\undefined
\ifx\selectfont\undefined
\documentstyle[11pt]{article}
\else
\documentstyle[oldlfont,11pt]{article}
\fi
\else
\documentclass[11pt]{article}
\usepackage{latexsym}
\fi
\font\tenmsb=msbm10  at 11pt
\font\sevenmsb=msbm7 at 9pt
\font\fivemsb=msbm5 at 7pt
\newfam\msbfam
\textfont\msbfam=\tenmsb
\scriptfont\msbfam=\sevenmsb
\scriptscriptfont\msbfam=\fivemsb
\def\Bbb#1{\fam\msbfam\relax#1}
        \def\nat{{\Bbb N}}
        \def\real{{\Bbb R}}
        \def\D{{\cal D}}
        \def\F{{\cal F}}
        \def\G{{\cal G}}
        \def\I{{\cal I}}
        \def\R{{\cal R}}
        \def\S{{\cal S}}
        \def\T{{\cal T}}
\def\supp{\mathop{\rm supp}\nolimits}
\def\ep{\varepsilon}
\def\Dstabs{$\Delta$-stabilizes}
\def\Dstab{$\Delta$-stabilize}
\def\gtrsim{\mathop{\buildrel >\over\sim}\nolimits}
\def\lesssim{\mathop{\buildrel <\over\sim}\nolimits}
\def\norm#1{{|\!|\!| #1 |\!|\!|}}
\def\ov{\overline}

\newtheorem{fact}{Fact}[section]
\newtheorem{thm}[fact]{Theorem}
\newtheorem{lemma}[fact]{Lemma}
\newtheorem{prop}[fact]{Proposition}
\newtheorem{cor}[fact]{Corollary}
\newtheorem{defin}[fact]{Definition}
\newtheorem{example}[fact]{Example}
\newtheorem{remark}[fact]{Remark}

\newtheorem{problem}[fact]{Problem}
\newenvironment{proof}{\medskip \par \noindent {\it Proof.}\ }{\hfill
$\Box$ \medskip \par}
\def\pf#1{\medskip \par \noindent {\it #1.}\ }
\def\endpf{\hfill $\Box$ \medskip \par}
\begin{document}
\title{Proximity to $\ell_1$ and Distortion in Asymptotic $\ell_1$
Spaces}
\author{Edward  Odell \thanks{Research partially supported by
NSF and TARP.}
\and Nicole Tomczak-Jaegermann \thanks{Research partially
  supported by NSERC.}
\and Roy  Wagner \thanks{Research partially   supported by BSF.}
}
\date{}
\maketitle
\begin{abstract}
  For an asymptotic $\ell_1$ space $X$ with a basis  $(x_i)$
  certain asymptotic $\ell_1$
  constants, $\delta_\alpha (X)$ are defined  for $\alpha <\omega_1$.
  $\delta_\alpha (X)$ measures the equivalence between all
  normalized block bases $(y_i)_{i=1}^k$ of $(x_i)$ which are
  $S_\alpha$-admissible with respect to $(x_i)$ ($S_\alpha$ is the
  $\alpha^{th}$-Schreier class of sets) and the unit vector basis of
  $\ell_1^k$. This leads to the concept of the delta spectrum of $X$,
  $\Delta (X)$, which  reflects the behavior of stabilized limits of
  $\delta_\alpha (X)$.  The analogues of these constants under all
  renormings of $X$ are also defined and studied.
  We investigate  $\Delta (X)$ both in general and for spaces
  of bounded distortion. We also prove several  results on
  distorting the classical Tsirelson's space $T$ and its relatives.
\end{abstract}

\baselineskip=18pt                      

\section{Introduction}

The first non-trivial example of what is now called an asymptotic
$\ell_1$ space was discovered by Tsirelson \cite{Ts}. This space and its
variations were extensively studied in many papers (see
\cite{CS}).  While the finite-dimensional asymptotic structure of
these spaces is the same as that of $\ell_1$, they do not contain an
infinite-dimensional subspace isomorphic to $\ell_1$, and thus their
geometry is inherently different.

The idea of investigating the geometry of a Banach space by studying
its asymptotic finite-dimensional subpaces arose naturally in recent
studies related to problems of distortion, i.e. the stabilization of
equivalent norms on infinite dimensional subspaces of a given Banach
space.  These ideas were further developed and precisely formulated
in \cite{MMT}.

By a finite-dimensional asymptotic subspace of $X$ we mean a subspace
spanned by blocks of a given basis living sufficiently far along the
basis.  By an asymptotic $\ell_p$ space we mean a space all of whose
asymptotic subspaces are $\ell_p^n$, i.e. any $n$ successive
normalized blocks of the basis $\{e_i\}_{i=1}^{\infty}$ supported
after $e_n$ are $C$-equivalent to the unit vector basis of $\ell_p^n$.

In this paper we introduce  a concept which bridges the gap
between this ``first order'' structure of an asymptotic $\ell_1$ space
and the global structure of its infinite-dimensional subspaces. This
concept employs a hierarchy of families of finite subsets of $\nat$
of increasing complexity, the Schreier classes $(\S_\alpha)_{\alpha <
{\omega_1}}$ introduced in \cite{AA}.  For $\alpha < \omega_1$ we
define what it means for a normalized block basis to be
$\S_\alpha$-admissible with respect to the basis $(e_i)$, and then
measure the equivalence constant between all such blocks and the
standard unit vector basis of $\ell_1$, obtaining the parameter
$\delta_\alpha(e_i)$. These constants increase when passing to block
bases and this leads us to define the $\Delta$-spectrum of $X$,
$\Delta(X)$, to be the set of all stabilized limits $\gamma =
(\gamma_\alpha)$ of $(\delta_\alpha(e_i))$ as $(e_i)$ ranges over all
block bases of $X$.

We show that these concepts provide useful and efficient tools for
studying the infinite dimensional and asymptotic structure of
asymptotic $\ell_1$ spaces.  Indeed, even some first order asymptotic
problems require a higher order analysis.  The behavior of the
$\Delta$-spectrum of $X$ has deep implications in regard to the
distortability of $X$ and its subspaces.

\medskip

We now describe the contents of the paper in more
detail.

\smallskip

Section 2 reviews concepts and results concerning distortion and
asymptotic $\ell_1$ spaces.  We sketch the proof of the
2-distortability of Tsirelson's space in Proposition~\ref{B.d}. This
leads to a natural question as to whether the asymptotic structure of
$T$ can be distorted: can $T$ be given an equivalent norm such that
its asymptotic subspaces are closer to $\ell_1$?  Without resorting to
the higher order analysis developed in subsequent sections we only
obtain a partial solution (the complete solution is then provided in
Section 5).

In Section 3 we define the Schreier families $\S_\alpha$ and establish
some facts about their mutual relationship which are crucial for  our
later work.

Section 4  contains precise definitions of all the asymptotic $\ell_1$
constants which we introduce in this  paper. We also define the spectrum
$\Delta(X)$. Elements $\gamma = (\gamma_\alpha)_{\alpha < \omega_1}$
of the spectrum satisfy $\gamma_\alpha \gamma_\beta \le \gamma_{\alpha
+ \beta}$ for all $\alpha$, $\beta < \omega_1$
(Proposition~\ref{D.j}). It follows that
$\hat\gamma_\alpha = \lim_{n \to \infty}\gamma_{\alpha \cdot n}^{1/n}$
exists for all $\alpha < \omega_1$  and it is shown to equal
$\ddot\delta_\alpha (Y)$  for some subspace $Y \subseteq X$
(Proposition~\ref{D.m}).
$\ddot\delta_\alpha (Y)$ is defined to be the largest of
$\delta_\alpha ((x_i), |\cdot|)$ as $(x_i)$ ranges over all block bases
of $Y$ and $|\cdot|$ over all equivalent norms. The constants
$(\ddot\delta_\alpha(X))_{\alpha < \omega_1}$ exhibit a remarkable
regularity. They are constantly one until $\alpha$ reaches
the spectral index of $X$, $I_\Delta (X)$; and then decrease
geometrically to 0 as $\alpha$ reaches $I_\Delta (X) \cdot \omega$
(Theorem~\ref{F.c}).
An important tool in this section is the  renorming result
of Theorem~\ref{D.nd}.

Section 5 contains the calculation of  asymptotic constants for
various asymptotic $\ell_1$ spaces. We consider $T$ along with various
other Tsirelson and mixed Tsirelson spaces. These and other examples
show that there is potentially considerable variety in the spectrum of
$X$ despite the regularity conditions imposed when considering all
renormings. In addition it is shown that for $\gamma \in \Delta (X)$
an appropriate block basis in $X$ admits a lower $T_{\gamma_1}$
block Tsirelson estimate.

The central theme of Section 6 is the following problem: Does there
exist an asymptotic $\ell_1$ Banach space  of bounded distortion?
In particular, is Tsirelson's space of bounded distortion?  We apply
our work to obtain some partial results in this and related directions.
We consider the consequences of assuming
that an asymptotic $\ell_1$ space is of bounded
distortion. In particular the asymptotic constants must
behave in a geometric  fashion (Theorem~\ref{F.f},
Corollary~\ref{F.fa}, Propositions~\ref{F.i} and \ref{F.j}).
Also, an asymptotic $\ell_1$ space of  bounded  distortion
bears a striking resemblance to a subspace of a Tsirelson-type
space $T(\S_{\alpha},\theta)$  for some 
$ \alpha < \omega_1 $  and $0 < \theta < 1$ (Theorem~\ref{F.l}).
Furthermore we show that a renorming of Tsirelson's space $T$
for which there exists $\gamma$ in the spectrum with $\gamma_1 = 1/2$
cannot distort $T$ by more than a fixed constant (Theorem~\ref{E.b}).

\section{Preliminaries}

In this paper we shall use certain  notation and basic facts from
Banach space theory, as presented in \cite{LT}.  Furthermore,
$X,Y,Z,\ldots$ shall denote separable infinite dimensional Banach
spaces.  By $Y\subseteq X$ we mean that $Y$ is a closed
infinite-dimensional linear subspace of $X$.
By $S(X) = \{x\in X: \|x\| =1\}$ we denote the unit sphere of $X$.

If $(e_i)$ is a basic sequence and $F\subseteq \nat$, $\langle
e_i\rangle_{ F}$ is the linear span of $\{e_i:i\in F\}$ and
$[e_i]_{ F}$ is the closure of $\langle e_i\rangle_{i\in F}$.  For
$F,G\subseteq \nat$ the notation  $F<G$ means  that $\max F<\min G$
or either $F$ or $G$ is empty.
$F<G$ are {\em adjacent intervals of\/} $\nat$ if for some
$k \le m < n$, $F = [k,m]= \{i \in \nat: k \le i \le m\}$
and $G = [m+1, n]$.
If $x\in \langle e_i\rangle$ and  $x= \sum a_i e_i$
then $\supp (x) = \{i:a_i\ne 0\}$ is the support of $x$ with respect
to $(e_i)$ (w.r.t. $(e_i)$).  For $x,y\in \langle
e_i\rangle$, we write $x<y$ if  $\supp (x) < \supp (y)$.  By
$(x_i)\prec (e_i)$ we shall mean that $(x_i)$ is a block basis of
$(e_i)$.  We say that $Y$ is a {\em block subspace\/} of $X$,
$Y\prec X$, if $X$ has a basis $(x_i)$ and
$Y= [y_i]_{ \nat}$ for some $(y_i)\prec (x_i)$.

\subsection{Distortion}
If a Banach space $(X,\|\cdot\|)$ is given an equivalent norm $|\cdot|$
we define the distortion of $|\cdot|$ by
\begin{defin}\label{A.a}
  $$d(X,|\cdot|) = \inf_Y\sup\left\{\frac{|x|}{|y|} :x,y\in
    S(Y,\|\cdot\|)\right\}\ ,$$
  where the infimum is taken over all infinite-dimensional subspaces
  $Y$ of $X$.
\end{defin}
\begin{remark}\label{A.ad}\rm
  If $X$ has a basis, then a standard approximation argument easily shows
  that in the above formula for $d(X,|\cdot|)$ it  is sufficient to take
  the infimum over all block subspaces $Y \prec X$; and this is the form
  of the definition we shall always use.
\end{remark}

The parameter $d(X,|\cdot|)$ measures how close $|\cdot|$  can be made
to being a multiple of $\|\cdot\|$, by  restricting  to an
infinite-dimensional subspace.

\begin{defin}\label{A.b}
  For $\lambda>1$, $(X,\|\cdot\|)$ is {\em $\lambda$-distortable\/} if
  there exists an equivalent norm $|\cdot|$ on $X$ so that
  $d(X,|\cdot|) \ge \lambda$.  $X$ is {\em distortable\/}
  if it is $\lambda$-distortable for some $\lambda >1$.  $X$ is {\em
    arbitrarily distortable\/} if it is $\lambda$-distortable for all
  $\lambda >1$.
\end{defin}

\begin{defin}\label{A.c}
  A space $(X,\|\cdot\|)$ is of {\em $D$-bounded distortion\/} if for all
  equivalent norms $|\cdot|$ on $X$ and all $Y\subseteq X$,
  $d(Y,|\cdot|)\le D$. A space $X$ is of {\em bounded distortion\/} if
  it is of $D$-bounded distortion for some $D < \infty$.
\end{defin}

Let us mention a more geometric approach to distortion. A subset $A
\subseteq X$ is called {\em asymptotic\/} if $\mbox{\rm dist}(A, Y)
=0$ for all infinite-dimensional subspaces $Y$ of $X$, i.e. for all
$Y$ and $\ep >0$ there is $x \in A$ such that $\inf_{y\in Y} \|x-y\| <
\ep$. Given $\eta >0$, consider the following property of $X$: there
exist $A, B \subseteq S(X)$ and $A^*$ in the unit ball of $X^*$ such
that: (i) $A$ and $B$ are asymptotic in $X$; (ii) for every $x\in A$
there is $x^* \in A^*$ such that $|x^* (x) | \ge 1/2$; (iii) for all
$y \in B$ and $x^* \in A^*$, $|x^* (y)| <\eta$. It is well known and
easy to see that if $d(X,|\cdot|)\ge \lambda$ for some equivalent norm
$|\cdot|$ on $X$ then in some $Y\subseteq X$ there exist such
asymptotic (in $Y$) ``almost biorthogonal'' sets, with $\eta =
1/\lambda$.  Conversely, given sets $A, B$ and $A^*$ as above, let
$|x| = \|x\| + (1/\eta)\sup \{|x^* (x)| : x^* \in A^*\}$ for $x \in
X$. Then $d(X,|\cdot|)\ge (1/2 + 1/4\eta)$.

\medskip
A  proof of the following simple  proposition is left for the reader.
Part b) was shown in \cite{To}.

\begin{prop}\label{A.d}
\begin{description}
\item[a)] Let $(X,\|\cdot\|)$ be of $D$-bounded distortion and let
  $|\cdot|$ be an equivalent norm on $X$. Then for all $\ep >0$ and
  $Y\subseteq X$ there exists $Z\subseteq Y$ and $c>0$ so that $|z|\le
  c\|z\| \le (D+\ep) |z|$ for all $z\in Z$.
\item[b)] Every Banach space  contains either an arbitrarily
  distortable subspace or a subspace of bounded distortion.
\end{description}
\end{prop}

Note that if $X$ has a basis then one may replace $Y\subseteq X$
and $Z\subseteq Y$,  in Definition~\ref{A.c}  and Proposition~\ref{A.d},
and the definition of an asymptotic set,  by
$Y\prec X$ and $Z\prec Y$, respectively.

It was shown in \cite{OS1}, \cite{OS2} that every $X$ contains either
a distortable subspace or a subspace isomorphic to $\ell_1$ or $c_0$
(both of which are not distortable \cite{Jam}).  Currently no examples
of distortable spaces of bounded distortion are known.  It is known
that such a space would for some $1\le p\le\infty$ necessarily contain
an asymptotic $\ell_p$ subspace (defined below for $p=1$) with an
unconditional basis and must contain $\ell_1^n$'s uniformly
(\cite{MiT}, \cite{Ma}, \cite{To}).

In light of these results it is natural to focus the search for a
distortable space of bounded distortion on asymptotic $\ell_1$ spaces
with an unconditional basis.

\subsection{Asymptotic $\ell_1$ Banach spaces}

Several definitions of asymptotic $\ell_1$ spaces appear in the
literature.  We shall use the definition  from \cite{MiT}.
\begin{defin}\label{B.a}
  A space $X$ with a basis $(e_i)$ is an {\em asymptotic $\ell_1$ space\/}
  (w.r.t. $(e_i)$) if there exists $C$ such that for all $n$ and all
  $e_n\le x_1<\cdots <x_n$,
  $$\| \sum_1^n x_i \| \ge (1/C) \sum_1^n \|x_i\|\ . $$
  The infimum of all $C$'s as above is called the {\em asymptotic $\ell_1$
  constant} of $X$.
\end{defin}

It should be noted that this definition depends on the choice of a basis:
a space $X$ may be asymptotic $\ell_1$ with respect to one basis but not
another.  However when the basis is understood, the reference to it is
often dropped.

In \cite{MMT} a notion of asymptotic structure of an arbitrary Banach
space was introduced; in as much as we shall not use it here, we omit
the details.  This led, in particular, to a more general concept of
asymptotic $\ell_1$ spaces; and spaces satisfying Definition~\ref{B.a}
above were called there ``stabilized asymptotic $\ell_1$''. Several
connections between the ``$MMT$-asymptotic structure'' of a space
\cite{MMT}  and
the ``stabilized asymptotic structure'' of its subspaces can be
proved; for instance,  an $MMT$-asymptotic $\ell_1$ space contains
an  asymptotic $\ell_1$ space in the sense of Definition~\ref{B.a}.

Before proceeding we shall briefly consider the prime example of an
asymptotic $\ell_1$ space not containing $\ell_1$, namely Tsirelson's
space $T$ \cite{Ts}.  Our discussion will motivate our subsequent
definitions.  The space $T$ is actually the dual of Tsirelson's
original space.  It was described in \cite{FJ} as follows.

Let $c_{00}$ be the linear space of finitely supported sequences.  $T$
is the completion of $(c_{00},\|\cdot\|)$ where $\|\cdot\|$ satisfies
the implicit equation
$$\|x\| = \max \left( \|x\|_\infty ,\sup \biggl\{ \frac12 \sum_{i=1}^n
  \|E_ix\| : n\in \nat \mbox{ and } n\le E_1 <\cdots <E_n
  \biggr\}\right)\ .$$

In this definition the $E_i$'s are finite subsets of $\nat$.  $E_ix$ is
the restriction of $x$ to the set $E_i$.  Thus if $x= (x(j))$ then
$E_ix(j) = x(j)$ if $j\in E_i$ and $0$ otherwise.  Of course it must be
proved that such a norm exists.  The unit vector basis $(e_i)$ forms a
1-unconditional basis for $T$ and $T$ is reflexive. 
If $e_n \le x_1<\cdots <x_n$ w.r.t.
$(e_i)$ then $\|\sum_1^n x_i\| \ge \frac12 \sum_1^n \|x_i\|$ and so $T$
is asymptotic $\ell_1$ with constant less than or equal to 2.  The next
proposition is the best that can currently be said about distorting $T$.
The proof, which we sketch, is illustrative.

\begin{prop}\label{B.d}
  $T$ is $(2-\ep)$-distortable for all $\ep>0$.
\end{prop}
\begin{proof}
  {\it (Sketch)\ \ } Let $\ep>0$ and choose $n$ so that $1/n <\ep$.  Define
  for $x\in T$,
  $$|x| = \sup \biggl\{ \sum_{i=1}^n \|E_ix\| : E_1
  <\cdots<E_n\biggr\}\ .$$
  Clearly, $\|x\| \le |x| \le n\|x\|$ for  $x \in T$ (in fact,
  for $n \le x$, $|x| \le 2\|x\|$).  Let $(x_i)\prec
  (e_i)_n^\infty$.  For any $k>n$ some normalized sequence $(y_i)_1^k
  \prec (x_i)_k^\infty$ is equivalent to the unit vector basis of
  $\ell_1^k$, with the equivalence constant as close to 1 as we wish.
  Thus if $y=(1/k) \sum_1^k y_i$, then $\|y\|\approx 1$.  Also if
  $E_1<\cdots < E_n$ then setting $I= \{i:E_j\cap \supp (y_i) \ne
  \emptyset$ for at most one $j\}$ and $J= \{1,\ldots,k\}\setminus I$ we
  have that $|J| \le n$ and
\begin{eqnarray*}
  \sum_1^n \|E_jy\| & \le & \frac1k \biggl( \sum_{i\in I} \|y_i\|
  + \sum_{i\in J} \sum_j\|E_j y_i\|\biggr)\\
  & \le & \frac1k \biggl( \sum_{i\in I} \|y_i\|
  + \sum_{i\in J} 2\|y_i\|\biggr)\\
  &\le & \frac1k (k-|J| + 2|J|) \le 1+\frac{n}k.
\end{eqnarray*}
Thus  $\inf \{|x| : \|x\|=1,\ x\in \langle x_i\rangle\} =1$.

Now let $z= (2/n) \sum_1^n z_i \in \langle x_i\rangle_n^\infty$
where $z_1 < \cdots < z_n$ and each $z_i$ is an $\ell_1^{k_i}$-average
of the sort just considered.  Here $k_{i+1}$ is taken very large
depending on $\max\supp (z_i)$ and $\ep$.  Since $\|z_i\|\approx 1$,
it follows that $|z| \ge (2/n) \sum \|z_i\| \approx 2$.  Yet, if
$m\le E_1<\cdots <E_m$, and $i_0$ is the smallest $i$ such that $e_m <
\max\supp(z_i)$, then the growth condition for $k_i$ implies that
$k_i$ is much larger than $m$ for $i_0 < i \le n$. Hence by the
argument above
\begin{eqnarray}\nonumber
  \frac12 \sum_1^m \|E_j z\| & = &   \frac1n
       \sum_1^m \biggl\|E_j \biggl(\sum_1^n z_i\biggr)\biggr\|\\
\label{B.form}
  &\le & \frac1n \biggl( \sum_{j=1}^m \|E_j z_{i_0}\|
  + \sum_{i={i_0}+1}^n \sum_{j=1}^m \|E_j z_i\|\biggr)\\
  &\le & \frac1n \biggl( 2\|z_{i_0}\|
  + \sum_{i={i_0}+1}^n (1 + n/k_i)\biggr)
   \lesssim  \frac{n+1}n < 1+\ep\ .\nonumber
\end{eqnarray}
By the definition of the norm we get  $\|z\| \le 1 +\ep$.
This implies  $\sup \{|z| : \|z\|=1,\ z\in \langle x_i\rangle\} >
2/(1+\ep)$.
\end{proof}

Later we shall say that a sequence $(y_i)_1^k$ is $\S_1$-admissible w.r.t.
$(x_i)$ if $x_k \le  y_1 <\cdots < y_k$.  In the above proof we needed to
consider an admissible sequence of admissible sequences; what we shall
later call $\S_2$-admissible.

Inequality (\ref{B.form}) obviously shows that the asymptotic $\ell_1$
constant of $T$ is greater  than or equal, and hence equal, to 2.
Furthermore, if $X\prec T$, then $X$ is an asymptotic $\ell_1$ space
with constant again equal to 2.  In other words, passing to a block
basis of $T$ does not improve the asymptotic $\ell_1$ constant.
Vitali Milman asked the question what would happen if in addition we
renormed?  The above technique gives that the constant cannot be
improved too much.

\begin{prop}\label{B.e}
  If $|\cdot|$ is any equivalent norm on $X\prec T$ then $X$ is
  asymptotic $\ell_1$ with constant at least $ \sqrt2$.
\end{prop}
\begin{proof}
  {\it (Sketch)\ \ }  Let  $X\prec T$ and consider an equivalent norm
  $|\cdot|$ on $X$ so that $(X,|\cdot|)$ is asymptotic $\ell_1$ with
  constant $\theta $.  By multiplying $|\cdot|$ by a constant
  and passing to a block subspace of $X$ if necessary we may assume
  that $\|\cdot\| \ge |\cdot|$ on $X$ and for all $Y\prec X$ there
  exists $y\in Y$ with $\|y\|=1$ and $|y|\approx 1$.  Given $n$,
  choose $z_1<z_2<\cdots <z_n$ w.r.t. $X$ so that $z_i =(1/k_i)
  \sum_1^{k_i} z_{i,j}$ where $z_{i,1}<\cdots < z_{i,k_i}$ in $X$ and
  $\|z_{i,j}\|=1\approx |z_{i,j}|$.  Here $k_{i+1}$ is again large
  depending upon $z_i$.
 
  Let $z= (2/n) \sum_1^n z_i$.  Then as before we obtain $|z| \le
  \|z\| \lesssim 1+ (1/n)$.  On the other hand, $|z| \ge (2/n \theta)
  \sum_1^n |z_i| \gtrsim (2/n \theta^{2})n = 2/\theta^{2}$.
  Hence $2/\theta^2 \lesssim 1$.
\end{proof}

\begin{remark}\label{B.ea}\rm
  For any $0 < \theta <1$ Tsirelson's space $T_\theta$ is defined  by
  the implicit equation analogous to the definition of $T$, in which
  the constant $1/2$ is replaced by $\theta$.   The properties of
  $T$ remain valid for $T_\theta$ as well, with appropriate  modification
  of the constants involved.
\end{remark}

These results indicate that it could be of advantage to consider the
$\ell_1$-ness of sequences which are $\S_2$-admissible with respect to
a basis or even $\S_n$-admissible.  We do so in this paper and we shall
obtain the best possible improvement of Proposition~\ref{B.e} in
Theorem~\ref{E.a} (see also Remark~\ref{E.ac}).  Of course the
beautiful examples of Argyros and Deliyanni \cite{AD} of arbitrarily
distortable mixed Tsirelson spaces (described below) also show the
need for consideration of such notions when studying asymptotic
$\ell_1$ spaces.  Our point here is that these are needed even to
answer $\S_1$-admissibility questions.

\section{The Schreier families  $\S_\alpha$}

Let $\F$ be a set  of finite subsets of $\nat$.  $\F$ is {\em
  hereditary\/} if whenever $G\subseteq F\in \F$ then $G\in \F$.  $\F$
is {\em spreading\/} if whenever $F= (n_1, \cdots, n_k)\in\F$,
with $n_1<\cdots <n_k$ and $m_1<\cdots<m_k$ satisfies
$m_i\ge n_i$ for $i\le k$ then $(m_1,\ldots,m_k)\in \F$.  $\F$ is {\em
  pointwise closed\/} if $\F$ is closed in the topology of pointwise
convergence in $2^{\nat}$.
A set $\F$ of finite subsets of $\nat$  having all three properties we call
{\em regular\/}.  If $\F$ and $\G$ are regular
we let
$$
\F[\G] =
\biggl\{\bigcup_1^n G_i: n\in\nat,\, G_1<\cdots < G_n,\, G_i \in \G
\mbox{\rm \ for } i\le n,\, (\min G_i)_1^n \in\F\biggr\}\ .
$$
Note that this operation satisfies the natural associativity condition
$(\F[\G_1])[\G_2] = \F [\G_0]$, where $\G_0 = \G_1[\G_2]$.

If $N= (n_1,n_2,\ldots)$ is a
subsequence of $\nat$ then $\F(N) = \{(n_i)_{i\in F} :F\in\F\}$. 
If  $\F$ is regular and $M$ is a subsequence of $N$ then,
since $\F$ is spreading, $\F(M) \subset \F(N)$.
If  $\F$ is regular and $n\in\nat$ we define $[\F]^n$ by
$[\F]^1 = \F$ and $[\F]^{n+1} = \F\Bigl[[\F]^n\Bigr]$. 
Finally, if $F$ is a finite
set,  $|F|$ denotes the cardinality of $F$.

\begin{defin}\label{C.a}
  {\rm \cite{AA}}
  The {\em Schreier classes\/} are defined by $\S_0 =
  \{\{n\}:n\in\nat\}\cup \{\emptyset\}$,
  $\S_1 = \{F\subseteq\nat :\min F\ge |F|\} \cup \{\emptyset\}$;
  for $\alpha<\omega_1$,
  $\S_{\alpha+1} = \S_1 [\S_\alpha]$, and if $\alpha$ is a limit ordinal we
  choose $\alpha_n\uparrow \alpha$ and set
  $$\S_\alpha = \{F:\mbox{\rm  for some } n\in\nat, \ F\in \S_{\alpha_n}
  \mbox{\rm  and } F\ge n \}\ .$$
\end{defin}

It should be noted that the definition of the $\S_\alpha$'s for
$\alpha\ge\omega$ depends upon the choices made at limit ordinals but
this particular choice is unimportant for our purposes.  Each
$\S_\alpha$ is a regular class of sets.  It is easy to see that
$\S_1\subseteq \S_2\subseteq \cdots$ and $\S_n[\S_m]= \S_{m+n}$, for
$n, m \in \nat$, but this fails for higher ordinals.  However we do
have

\begin{prop}\label{C.b}
\begin{description}
\item[a)] Let $\alpha<\beta<\omega_1$. Then there exists $n\in\nat$ so
  that if $n\le F\in \S_\alpha$ then $F\in \S_\beta$.
\item[b)] For all $\alpha,\beta<\omega_1$ there exists a subsequence
  $N$ of $\nat$ so that
  $\S_\alpha [\S_\beta] (N) \subseteq \S_{\beta+\alpha}$.
\item[c)] For all $\alpha,\beta<\omega_1$ there exists a subsequence
  $M$ of $\nat$ so that
  $\S_{\beta+\alpha} (M) \subseteq \S_\alpha [\S_\beta]$.
\end{description}
\end{prop}
We start with  an easy formal observation.
\begin{lemma}\label{C.ba}
  Let $\F$and $\G$ be sets of finite subsets of $\nat$ and let $\G$ be
  spreading. Assume that there exists a subsequence $N$ of $\nat$ so that
  $\F(N) \subseteq \G$. Then for all  subsequences $L$ of $\nat$ there
  exists a subsequence $L'$ of $L$ with $\F(L') \subseteq \G$.
\end{lemma}
\begin{proof}
  Let $L = (l_i)$. Let $N = (n_i)$ such that $\F(N) \subseteq \G$.
  Since $\G$ is spreading, any $L' = (l_i') \subseteq (l_i)$ such that
  $l_i' \ge n_i$ for all $i$ satisfies the conclusion (for instance
  one can take $L'= (l_{n_i})$).
\end{proof}
\medskip \pf{Proof of  Proposition~\ref{C.b}}
  a) We proceed by induction on $\beta$.  If $\beta = \gamma+1$ then
  $\alpha\le\gamma$ and so we may choose $n$ so that if $n\le F\in
  \S_\alpha$ then $F\in \S_\gamma \subseteq \S_\beta$.  If $\beta$ is a
  limit ordinal and $\beta_n\uparrow \beta$ is the sequence used in
  defining $\S_\beta$, choose $n_0$ so that $\alpha<\beta_{n_0}$.  Choose
  $n\ge n_0$ so that if $n\le F\in \S_\alpha$ then $F\in \S_{\beta_{n_0}}$.
  Thus also $F\in \S_\beta$.
 
  b) We induct on $\alpha$.  Since $\S_0 [\S_\beta] = \S_\beta$, the
  assertion is clear for $\alpha=0$.  If $\alpha = \gamma+1$, then
  $\S_\alpha [\S_\beta] = \S_1[\S_\gamma [\S_\beta]]$ and $\S_{\beta+\alpha}
  = \S_1 [\S_{\beta+\gamma}]$.  Thus we can take $N$ to satisfy
  $\S_\gamma [\S_\beta] (N)\subseteq \S_{\beta+\gamma}$.
 
  If $\alpha$ is a limit ordinal we argue as follows. First, by
  Lemma~\ref{C.ba}, the inductive hypothesis implies that for
  every $\alpha' < \alpha$ and  every
  subsequence $L$ of $\nat$ there exists a subsequence $N$ of $L$ with
  $\S_{\alpha'}[\S_\beta] (N)\subseteq \S_{\beta+\alpha'}$.  Let
  $\alpha_n\uparrow \alpha$ and $\gamma_n\uparrow \beta+\alpha$ be the
  sequences of ordinals used to define $\S_\alpha$ and
  $\S_{\beta+\alpha}$, respectively.
 
  Choose subsequences of $\nat$, $L_1\supseteq L_2\supseteq \cdots$ so
  that $\S_{\alpha_k} [\S_\beta] (L_k) \subseteq \S_{\beta+\alpha_k}$.
  If $L_k = (\ell_i^k)_{i=1}^\infty$ we let $L$ be the diagonal
  $L=(\ell_k) = (\ell_k^k)_{k=1}^\infty$.  It follows that if $F\in
  \S_{\alpha_k} [\S_\beta] (L)$ and $F\ge \ell_k$ then $F\in
  \S_{\alpha_k} [\S_\beta] (L_k)$ and so $F\in \S_{\beta+\alpha_k}$.
  For each $k$ choose $\bar n(k)$ so that $\beta+\alpha_k
  <\gamma_{\bar n(k)}$.  Using a) choose $j(1) < j(2) <\cdots$ so that
  if $j(k)\le F\in \S_{\beta+\alpha_k}$ then $F\in \S_{\gamma_{\bar
  n(k)}}$.  Let $N= (n(k))_{k=1}^\infty$ be a subsequence of $\nat$
  with $n(k) \ge \ell_k \vee j(k) \vee \bar n(k)$ for all $k$. 
  Then if $F\in \S_\alpha [\S_\beta] (N)$ there exists $k$ so that
  $n_k\le F \in \S_{\alpha_k} [\S_\beta] (N)$ and so $ \ell_k \le F
  \in \S_{\beta+\alpha_k}$ and $j(k) \le F\in \S_{\gamma_{\bar n(k)}}$,
  whence  since $\bar n (k) \le F$  we have   $F\in \S_{\beta+\alpha}$.

  c) As in b) we induct on $\alpha$.  The cases $\alpha=0$   and
  $\alpha=\gamma+1$ are trivial.  Thus assume that $\alpha$ is a limit
  ordinal.  Let $\alpha_r \uparrow \alpha$ and $\bar\gamma_r\uparrow
  \beta+\alpha$ be the sequences defining $\S_\alpha$ and
  $\S_{\beta+\alpha}$ respectively.  We may write $(\bar\gamma_r) =
  (\bar\gamma_1,\ldots,\bar\gamma_{n_0-1},\beta+\gamma_{n_0},
  \beta+\gamma_{n_0+1},\ldots)$ where $\bar\gamma_i <\beta$ if
  $i<n_0$.  By a) there exists $m_0$ so that if $m_0 \le F\in
  \bigcup_1^{n_0-1} \S_{\bar\gamma_i}$ then $F\in \S_\beta$.  We shall
  take later $M= (m_i)_1^\infty$ where $m_1\ge m_0$.  By the inductive
  hypothesis and Lemma~\ref{C.ba}. choose sequences
  $L_{n_0} \supseteq L_{n_0+1} \supseteq \cdots$ so that
  $\S_{\beta+\gamma_k} (L_k) \subseteq \S_{\gamma_k} [\S_\beta]$ for
  $k\ge n_0$ and  $m_0\le L_{n_0}$.  If $L_k = (\ell_i^k)_i$ set $L=
  (\ell_k)$ where $\ell_k = \ell_k^k$ for $k\ge n_0$, and $m_0 \le
  \ell_1 <\cdots <\ell_{k-1} < \ell_k <\cdots$.  Thus if $k\ge n_0$,
  $\ell_k \le F\in \S_{\beta+\gamma_k}(L)$ implies that $F\in
  \S_{\gamma_k} [\S_\beta]$.  Also for $k<n_0$, $\ell_k\le F\in
  \S_{\bar\gamma_k}$ implies that $F\in \S_\beta$.  For $k\ge n_0$
  choose $\bar m(k)$ so that $\beta+\gamma_k <\beta + \alpha_{\bar
    m(k)}$.  By a) there exists $n(k)$ so that $n(k) \le F\in
  \S_{\gamma_k}[\S_\beta]$ implies that $F\in \S_{\alpha_{\bar
      m(k)}}[\S_\beta]$ for all $k\ge n_0$.  Finally we choose $M=
  (m(k))$ where $m(k) = \ell_k$ for $k<n_0$ and $m(k) \ge \ell_k
  \vee \bar m (k) \vee n(k)$ for $k\ge n_0$.  Thus if $F\in
  \S_{\beta+\alpha} (M)$ then $F\in \S_\beta$ or else there exists $k\ge
  n_0$ with $m(k) \le F\in \S_{\beta +\gamma_k}(M)$.  Hence $\ell_k \le F\in
  \S_{\gamma_k} [\S_\beta]$ and so $n(k) \le F\in \S_{\alpha_{\bar m(k)}}
  [\S_\beta]$.  Since $F\ge \bar m(k)$ we get that $F\in
  \S_\alpha[\S_\beta]$.
\endpf
\begin{cor}\label{C.d}
  For all $\alpha <\omega_1$ and $n\in\nat$ there exist subsequences
  $M$ and $N$ of $\nat$ satisfying $[\S_\alpha]^n(N) \subseteq
  \S_{\alpha\cdot n}$ and $\S_{\alpha\cdot n} (M) \subseteq
  [\S_\alpha]^n$.
\end{cor}
\begin{proof}
  This is easily established by induction on $n$ using
  Proposition~\ref{C.b}. For example, if $[\S_\alpha]^n(P) \subseteq
  \S_{\alpha\cdot n}$ and $\S_\alpha[\S_{\alpha\cdot n}](L) \subseteq
  \S_{\alpha\cdot (n+1)}$, let $N = (p_{l_i})$ (here $P = (p_i)$ and
  $L = (l_i)$). Then $[\S_\alpha]^n(N) \subseteq \S_{\alpha\cdot
    n}(L)$ and so $[\S_\alpha]^{n+1}(N) = \S_\alpha [\S_\alpha]^n(N)
  \subseteq \S_\alpha[\S_{\alpha\cdot n}](L) \subseteq
  \S_{\alpha\cdot(n+1)}$.
\end{proof}

\begin{remark}\label{C.e}\rm
  The Schreier family $\S_\alpha$ has been used in {\rm \cite{AA}} to
  construct an interesting subspace $S_\alpha$ of
  $C(\omega^{\omega^\alpha})$ as follows.  $S_\alpha$ is the completion
  of $c_{00}$ under the norm
  $$\|x\| = \sup \{ |\sum_{i \in E}x(i)| : E\in \S_\alpha\}\ .$$
  The unit vector basis is an  unconditional basis for $S_\alpha$.
  The space $S_\alpha$ does not embed into $C(\omega^{\omega^\beta})$ for
  any $\beta<\alpha$.
\end{remark}

The next important proposition is a slight generalization of a result
in \cite{AD}  and is a descendent of results in \cite{Belle}.

\begin{prop}\label{C.l}
  Let $\beta <\alpha< \omega_1$, $\ep >0$ and let $M$ be a subsequence
  of $\nat$.  Then there exists a finite set $F\subseteq M$ and
  $(a_j)_{j\in F}\subseteq \real^+$ so that $F\in \S_\alpha (M)$,
  $\sum_{j\in F} a_j=1$ and if $G\subseteq F$ with $G\in \S_\beta$ then
  $\sum_{j\in G} a_j<\ep$.
\end{prop}
\begin{proof}
  We proceed by induction on $\alpha$.  The result is clear for
  $\alpha=1$.  Let $M= (m_i)$.  We choose $1/k <\ep$, $F\subseteq M$,
  $F>m_k$, $|F|=k$ and let $a_j =1/k$ if $j\in F$.
 
  If $\alpha$ is a limit ordinal let $\alpha_n\uparrow \alpha$ be the
  sequence used to define $\S_\alpha$.  Choose $n$ so that $\beta
  <\alpha_n$.  Applying the induction hypothesis to $\beta$,
  $\alpha_n$ and $\{m\in M: m\ge m_n\}$ yields the result.
 
  If $\alpha = \gamma+1$ we may assume (by Proposition~\ref{C.b}) that
  $\beta=\gamma$.  If $\gamma$ is a limit ordinal let
  $\gamma_n\uparrow\gamma$ be the sequence used to define $\S_\gamma$.
  Choose $k$ so that $1/k <\ep/2$.  Choose sets $F_i\subseteq M$ with
  $m_k\le F_1<\cdots <F_k$ along with scalars $(a_j)_{j\in \bigcup_1^k
    F_i} \subseteq \real^+$ and $n_1<\cdots <n_k$ satisfying the
  following:
\begin{description}
\item[1)] $\sum_{j\in F_i} a_j =1$ for $i\le k$
\item[2)] $F_i \in \S_{\gamma_{n_i}} (M)$ and $m_{n_i} <  F_i$
  for $1 \le i\le k$
\item[3)] $\sum_{j\in G} a_j <1/2^i$ if $G\subseteq F_{i+1}$ with
  $G\in \S_{\gamma_\ell}$ whenever $\ell \le \max F_i$ for $1\le i<k$.
\end{description}

Let $F= \bigcup_1^k F_i$.  Then $F\subseteq M$ and $F\in
\S_{\gamma+1}(M)$.  For $j\in F$ set $b_j = k^{-1} a_j$.  Then
$(b_j)_{j\in F} \subseteq \real^+$, $\sum_{j\in F} b_j=1$ and if $G\in
\S_\gamma$, $G\subseteq F$ then $\sum_{j\in G} b_j<\ep$.  Indeed there
exists $n$ with $n\le G\in \S_{\gamma_n}$.  Thus if $i_0 = \min \{i:
G\cap F_i\ne\emptyset\}$ then $n \le \max F_{i_0}$ and so by 3),
$$\sum_{j\in G} b_j = \sum_{i=i_0}^k \biggl( \sum_{j\in F_i\cap G}
b_j\biggr) \le \frac1k \left( 1+ \frac1{2^{i_0}} +\cdots +
  \frac1{2^k}\right) <\ep\ .$$

If $\gamma = \eta+1$ we again choose $1/k <{\ep}/2$ and sets $m_k \le
F_1<\cdots <F_k$, $F_i \in \S_\gamma (M)$, along with $(a_j)_{j\in
  F_i} \subseteq \real^+$, $\sum_{j\in F_i} a_j=1$ so that if $G\in
\S_\eta$ then $\sum_{j\in G\cap F_{i+1}} a_j < (1/2^i \max (F_i))$ for
$1 \le i <k$.  As above we set $F=\bigcup_1^k F_i$ and let $b_j =
k^{-1}a_j$ if $j\in F_i$.  Thus $F\in \S_\alpha (M)$ and if $G\in
\S_\gamma$, $G\subseteq F$, write $G=\bigcup_{s=1}^p G_s$ where $p\le
G_1<\cdots <G_p$ and $G_i \in \S_\eta$ for each $i$.  Then if $i_0 =
\min \{i:G\cap F_i\ne\emptyset\}$, since $\max (F_{i_0})\ge p$,
$$\sum_{j\in G} b_j \le \frac1k + \sum_{i=i_0+1}^k \frac1{2^i}\
\frac{p}k\ \frac1{\max (F_{i})}\le \frac1k +  \frac1k \sum_{i=i_0+1}^k
2^{-i} < \ep \ .$$
\end{proof}

\begin{defin}\label{C.m}
  Let $\ep>0$ and $\beta<\alpha <\omega_1$.  If $(e_i)$ is a
  normalized basic sequence, $M$ is a subsequence of $\nat$ and $F$
  and $(a_i)_{i\in F}$ are as in Proposition~\ref{C.l}, we call $x=
  \sum_{i\in F} a_i e_i$ an {\em $(\alpha,\beta,\ep)$-average of}
  $(e_i)_{i\in M}$.  If $(x_i)$ is a normalized block basis of $(e_i)$
  and $F$ and $(a_i)_{i\in F}$ are as in Proposition~\ref{C.l} for $M
  = (\min \supp (x_i))$, we call $x= \sum_{i\in F} a_i x_i$ an {\em
    $(\alpha,\beta,\ep)$-average of $(x_i)$ w.r.t. $(e_i)$}.
\end{defin}

The Schreier families are large within the set of all classes of
pointwise closed subsets of $[\nat]^{<\omega}$. Our next two
propositions show that they are in a sense the largest among all
regular classes of a given complexity.  To make this concept precise
we consider the index $I(\F)$ defined as follows.  Let $D(\F) = \{F\in
\F:$ there exist $(F_n)\subseteq \F$ with $1_{F_n}\to 1_F$ pointwise
and $F_n\ne F$ for all $n\}$, $D^{\alpha+1} (\F) = D(D^\alpha (\F))$
and $D^\alpha (\F) =\bigcap_{\beta<\alpha} D^\beta (\F)$ when $\alpha$
is a limit ordinal. Then
$$I(\F) = \inf \{\alpha <\omega_1 :D^\alpha (\F) = \{\emptyset\}\}\ .$$
$\F$ is a countable compact metric space in the topology of pointwise 
convergence and so $I(\F)$ must be countable, see e.g. \cite{Ku},
p.~261-262.

\begin{remark}\label{C.c}\rm
  The Cantor-Bendixson index of $\F$ (under the topology of pointwise
  convergence) is $I(\F)+1$.  This is because $\emptyset$ corresponds
  to the $0$ function and one needs one more derivative to get
  $\emptyset : D^{I(\F)+1} (\F) = \emptyset$, which defines the
  Cantor-Bendixson index.
\end{remark}

Now we have (\cite{AA})
\begin{prop}\label{C.n}
  For $\alpha <\omega_1$, $I(\S_\alpha) = \omega^\alpha$.
\end{prop}
\begin{proof}
  We induct on $\alpha$.  The result is clear for $\alpha=0$.  If the
  proposition holds for $\alpha$ it can be easily seen that for
  $n\in\nat$, $D^{\omega^\alpha\cdot n} (\S_{\alpha+1}) = \{F:$ there
  exists $k\in \nat$, $k>n$, with $F= \bigcup_1^{k-n} F_i$, $k\le F_1
  <\cdots < F_{k-n}$, and $F_i\in \S_\alpha $ for $i\le k-n\}$.  Hence
  $I(\S_{\alpha+1}) = \omega^{\alpha+1}$.  The case where $\alpha$ is a
  limit ordinal is also easily handled.
\end{proof}

\begin{prop}\label{C.j}
  If $\F$ is a regular set  of finite subsets of $\nat$ with $I(\F)\le
  \omega^\alpha$ then there exists a subsequence $M$ of $\nat$ with
  $\F(M)\subseteq \S_\alpha$.
\end{prop}

This proposition is a special case of  more complicated statements
(Proposition~\ref{C.k} and Remark~\ref{C.i}) below. First let us recall
(see e.g., \cite{M}) that every ordinal $\beta<\omega_1$ can be uniquely
written in Cantor normal form as
$$\beta = \omega^{\alpha_1}\cdot n_1 + \omega^{\alpha_2} \cdot n_2 +
\cdots + \omega^{\alpha_j}\cdot n_j$$
where $(n_i)_1^j \subseteq \nat$ and $\omega_1 >\alpha_1 >\cdots
>\alpha_j \ge0$.


\begin{defin}\label{C.g}
If $(\alpha_i)_1^j$ are countable ordinals and $(n_i)_1^j \subseteq
\nat$,  by $((\S_{\alpha_1})^{n_1},\ldots,(\S_{\alpha_j})^{n_j})$ we
denote  the class of subsets of $\nat$ that can be written in the form
$$E_1^1 \cup \cdots \cup E_{n_1}^1 \cup E_1^2 \cup \cdots\cup
E_{n_2}^2 \cup \cdots\cup E_1^j \cup\cdots\cup E_{n_j}^j$$ where
$E_1^1 < E_2^1 <\cdots < E_{n_j}^j$ and $E_i^k \in \S_{\alpha_k}$ for all
$i\le n_k$ and $k\le j$.
\end{defin}

\begin{prop}\label{C.k}
Let $\F$ be a regular set  of finite  subsets of $\nat$ with
$$I(\F) = \omega^{\alpha_1}\cdot n_1 + \cdots + \omega^{\alpha_{k-1}}
\cdot n_{k-1} + \omega^{\alpha_k}\cdot n_k,$$ in Cantor normal form.
Then there exists a subsequence $M$ of $\nat$ so that \newline
$\F(M) \subseteq ((\S_{\alpha_k})^{n_k},\ldots, (\S_{\alpha_1})^{n_1})$.
\end{prop}
\begin{remark}\label{C.i}\rm
The conclusion of the proposition holds even if $I(\F) <
\omega^{\alpha_1} \cdot  n_1 + \cdots + \omega^{\alpha_k}\cdot n_k$.
Indeed
this follows from the fact that if $\alpha<\beta$ and $\F(N) \subseteq
\S_\alpha$, and $N=(n_i)$, then there exists $r\in\nat$ so that
$\F((n_i)_{i\ge r}) \subseteq \S_\beta$ (by Proposition~\ref{C.b}(a)).
\end{remark}
\medskip \pf{Proof of  Proposition~\ref{C.k}}
We induct on $I(\F)$.  If $I(\F) =1$ then $\F$ contains only
singletons $\{n\}$ and so $\F(\nat) \subseteq ((\S_0)^1)$.

Assume the proposition holds for all classes with index $<\beta$,
and let $I(\F) = \beta$.  For $j\in \nat$ set $\F_j = \{F\in \F
:\{j\} \cup F\in \F$ and $j<F\}$.  Each $\F_j$ is regular.
					       
\noindent {\bf Case 1.} $\beta$ is a successor
Let $\beta = \omega^{\alpha_1} \cdot n_1 +\cdots + \omega^{\alpha_{k-1}}
\cdot n_{k-1} + n_k$    with $n_k>0$ ($\alpha_k=0$ here).
For every $j$,  $I(\F_j) \le\beta-1$.
Thus there exists $N_j\subseteq \nat$ such that
\begin{description}
\item[(1)] $\F_j(N_j) \subseteq ((\S_0)^{n_k-1},\ldots, (\S_{\alpha_1})^{n_1})\ ,$
\end{description}
with the  convention that  $(\S_0)^0 =\emptyset$.
Let $N_j=(n_i^j)_{i=1}^\infty$ and choose $N= (m_j)$ so that
$m_j\ge n_j^1 \vee \cdots \vee n_j^j$ for all $j$.

We shall show that
$$\F(N) \subseteq ((\S_0)^{n_k},\ldots,(\S_{\alpha_1})^{n_1})\ .$$

Indeed,  let $F\in \F$ and let  $\min F=j$ and  $G= F\setminus \{j\}$.
Then
$$(m_i)_{i\in F} = \{m_j\} \cup \{m_i\}_{i\in G} \in ((\S_0),(\S_0)^{n_k-1},
\ldots, (\S_{\alpha_1})^{n_1})\ ,$$
by (1), the choice of $N$ and
the fact that each $\S_\alpha$ is spreading.

\noindent {\bf Case 2.} $\beta$ is a limit ordinal.

Let $\beta = \omega^{\alpha_1}\cdot n_1 +\cdots + \omega^{\alpha_k}\cdot n_k$.
Note that $\alpha_k>0$.
We have $I(\F_j)<\beta$ for all $j$.

\noindent {\bf Case 2.1.} $\alpha_k$ is a successor.

Pick  $p_j\uparrow \infty$  such that for every $j$,
$$I(\F_j) \le \omega^{\alpha_1}\cdot n_1+\cdots + \omega^{\alpha_k}\cdot
(n_k-1) + \omega^{\alpha_k-1}\cdot p_j\ .$$

By induction there exist subsequences $N_j$ with
$$\F_j (N_j) \subseteq ((\S_{\alpha_k-1})^{p_j},(\S_{\alpha_k})^{n_k-1},\ldots,
(\S_{\alpha_1})^{n_1})\ .$$
Let $N_j = (n_i^j)$ and $N = (m_j)$ where $m_j \ge n_j^1 \vee \cdots \vee
n_j^j\vee (p_j +1)$ for all $j$.
If $F\in \F$ with $\min F=j$ and $G= F\setminus \{j\}$ then
$$(m_i)_{i\in G} = H_1 \cup \cdots \cup H_{p_j}\cup H$$
where $H_1<\cdots < H_{p_j} <H$; and $H_1,\ldots,H_{p_j}\in \S_{\alpha_k-1}$,
and $H\in ((\S_{\alpha_k})^{n_k-1},\ldots,  (\S_{\alpha_1})^{n_1})$.

Since $\{n_j\}\in \S_0\subseteq \S_{\alpha_k-1}$ and $m_j\ge p_j+1$,
we have, 
$\{m_j\} \cup H_1\cup \cdots \cup H_{p_j} \in \S_{\alpha_k}$. 
Thus
$$(m_i)_{i\in F} = \{m_j\} \cup H_1\cup\cdots\cup H_{p_j}\cup H
\in ((\S_{\alpha_k}),(\S_{\alpha_k})^{n_k-1},\ldots, (\S_{\alpha_1})^{n_1})\ .$$

\noindent {\bf Case 2.2.} $\alpha_k$ is a limit ordinal.

Let $\gamma_\ell\uparrow \alpha_k$ be the sequence of ordinals
defining $\S_{\alpha_k}$.
Set
$$\eta_\ell = \omega^{\alpha_1}\cdot n_1+\cdots + \omega^{\alpha_k}\cdot
(n_k-1) +\omega^{\gamma_\ell}\cdot 1$$
so that $\eta_\ell \uparrow \beta$.
Choose $\ell_j\uparrow \infty$ so that $I(\F_j) <\eta_{\ell_j}$.
As above choose $N_j$ so that $\F_j(N_j)\subseteq ((\S_{\gamma(\ell_j)})^1,
(\S_{\alpha_k})^{n_k-1},\ldots, (\S_{\alpha_1})^{n_1}$.
By Proposition~\ref{C.b}(a) there exists $r_j\in \nat$ so that
$r_j\le H\in \S_{\gamma(\ell_j)+1}$ implies $H\in \S_{\gamma(\ell_{j+1})}$.
Set $N_j = (n_i^j)$ and choose $N= (m_j)$ with
$m_j \ge n^1_j \vee \cdots \vee n_j^j\vee r_j \vee \ell_{j+1}$ for all $j$.

If $F\in \F$ with $\min F=j$ and  $G= F\setminus \{j\}$ then
$(m_i)_{i\in G} =  H_1 \cup H_2$  where $H_1 < H_2$, 
$H_1 \in \S_{\gamma(\ell_j)}$ and
$H_2 \in ((\S_{\alpha_k})^{n_k-1},\ldots, (\S_{\alpha_1})^{n_1})$.
Now $\{m_j\}\cup H_1 \in \S_{\gamma (\ell_j)+1}$ and by $m_j\ge r_j$
we have $\{m_j\} \cup H_1\in \S_{\gamma(\ell_{j+1})}$.
Also $m_j\ge \ell_{j+1}$ so $\{m_j\} \cup H_1\in \S_{\alpha_k}$.
\endpf

\begin{remark}\label{three.fourteen} \rm
The proof of Proposition~\ref{C.k} is due to Denny Leung and Wee Kee Tang.
They pointed out that our original proof was nonsense and supplied the
argument given.
We thank them for permission to reproduce it here.

In addition Denny Leung \cite{Leung}
has independently discovered a heirarchy of
sets similar to that of the Schreier classes.
\end{remark}

\begin{cor}\label{C.ka}
Let $\F$ be a pointwise closed class of finite subsets of $\nat$.
Then there exist $\alpha<\omega_1$ and a subsequence $M$ of $\nat$
so that  $\F(M) \subseteq \S_\alpha$.
\end{cor}
\begin{proof}
Let $\R$ be the regular hull of $\F$; that is, $\R = \{G: $ there
exists $F= (n_1,\ldots,n_k)\in \F$ with $G\subseteq (m_i)_1^k$ for some
$m_1<\cdots < m_k$ with $m_i \ge n_i$ for $i\le k\}$.

Clearly, $\R$ is hereditary and spreading.  We check that it is also
pointwise closed, and hence the corollary follows from
Proposition~\ref{C.j}.  Let $G_n\to G$ pointwise for some $(G_n)
\subseteq \R$.  If $|G| <\infty$ then $G$ is an initial segment of
$G_n$ for large $n$ and so $G\in \R$.  It remains to note that $|G|
=\infty$ is impossible.  If $G= (n_1,n_2,\ldots)$ then for all $k$,
$(n_1,\ldots,n_k)$ is a subset of some spreading of some set $F_k\in
\F$.  In particular $|\{ n\in F_k : n\le n_j\}| \ge j$ for $1\le j\le
k$.  Thus any limit point of $(F_k)_{k=1}^\infty$ is infinite which
contradicts the hypotheses that $\F$ is pointwise closed and consists
of finite sets.
\end{proof}

\begin{remark}\label{three.sixteen}\rm
R.~Judd \cite{Judd} has recently proved the following dichotomy result
for Schreier sets.
\end{remark}

\begin{thm}\label{rjudd}
Let $\F$ be a hereditary family of subsets of $\nat$ and let $\alpha <
\omega_1$.
Then either there exists a subsequence $M$ of $\nat$ so that $\S_\alpha (M)
\subseteq \F$ or there exist subsequences $M$ and $N$ of $\nat$ so that
$\F[M] (N) \subseteq \S_\alpha$,
where $\F [M] = \{F\subset M: F\in \F\}$. 
\end{thm}

For some other interesting properties of the Schreier classes we refer
the reader to \cite{AMT} and \cite{AO}.

\section{Asymptotic constants and $\Delta (X)$}

Asymptotic constants  considered in this paper will be determined
by the Schreier families  $\S_\alpha$; nevertheless it should be noted
that they can be introduced for a very general  class of families
of finite subsets of $\nat$.
\begin{defin}\label{D.a}
  If $\F$ is a regular set of finite subsets of $\nat$, a sequence of
  sets $E_1<\cdots < E_k$ is {\em $\F$-admissible\/} if $(\min
  (E_i))_{i=1}^k\in\F$.  If $(x_i)$ is a basic sequence in a Banach
  space and $(y_i)_1^k\prec (x_i)$, then $(y_i)_1^k$ is {\em
    $\F$-admissible\/} (w.r.t.  $(x_i)$) if $(\supp (y_i))_1^k$ is
  $\F$-admissible, where $\supp (y_i)$ is taken w.r.t. $(x_i)$.  We
  use a short form {\em $\alpha$-admissible} to mean
  $\S_\alpha$-admissible.
\end{defin}

The next definition was first introduced in \cite{Tom} for asymptotic
$\ell_p$ spaces with $1 \le p < \infty$.
\begin{defin}\label{D.c}
  Let $\F$ be a regular set  of finite subsets of $\nat$. For a basic
  sequence $(x_i)$ in a Banach space $X$ we define
  $\delta_{\F}(x_i)$ to be the supremum of $\delta\ge 0 $ such that
  whenever $(y_i)_1^k \prec (x_i)$ is  $\F$-admissible w.r.t. $(x_i)$
  then
  $$\|\sum_{i=1}^k y_i\| \ge \delta \sum_{i=1}^k \|y_i\|\ .$$ If $X$
  is a Banach space with a basis $(e_i)$ we write $\delta_{\F}(X)$ for
  $ \delta_{\F}(e_i)$.  For $\alpha <\omega_1$, we set $\delta_\alpha
  (x_i) = \delta_{\S_\alpha} (x_i)$ and $\delta_\alpha (X) =
  \delta_{\S_\alpha} (X)$.
\end{defin}

\begin{remark}\label{D.ca}\rm
  Note that   $\delta_{\F}(x_i)$ is equal to the supremum of all
  $\delta' \ge 0$  such that $\|y\| \ge \delta' \sum \|E_i y\|$,
  for all $y \in \langle x_i \rangle$ and all adjacent
  $\F$-admissible intervals    $E_1 < \cdots < E_k $ such that
  $\bigcup E_i \supseteq \supp (y)$. Here the support of $y$ and
  restrictions $E_i y$ are understood to be w.r.t. $(x_i)$.
  Indeed, clearly $\sup \delta' \ge \delta_{\F}(x_i)$. Conversely,
  given $(y_i)_1^k \prec (x_i)$   $\F$-admissible  we set $y = \sum
  y_i$  and we let $(E_1, \ldots, E_k)$  be adjacent intervals such
  that $E_i \supseteq \supp(y_i)$ and $\min E_i = \min \supp (y_i)$
  for all $i$.
\end{remark}

In as much as distortion problems involve passing to block subspaces
and renormings, it is natural to make two more definitions.

\begin{defin}\label{D.d}
  Let $\F$ be a regular set  of finite subsets of $\nat$ and let $(e_i)$
  be a basis for $X$.
\begin{eqnarray*}
  \dot\delta_{\F} (X) &=&\dot\delta_{\F}(e_i) = \sup
  \{\delta_{\F}(x_i)
  : (x_i) \prec (e_i)\} \quad   \mbox{and} \\
  \ddot\delta_{\F} (X) & = & \ddot\delta_{\F}(e_i) = \sup
  \{\dot\delta_{\F} ((e_i) ,|\cdot|) : |\cdot|\mbox{\rm\  is an equivalent
    norm on } X\}\ .
\end{eqnarray*}
We write $\dot\delta_{\S_\alpha} (X) = \dot\delta_\alpha(X)$ and
$\ddot\delta_{\S_\alpha} (X) = \ddot\delta_\alpha (X)$.
\end{defin}

The asymptotic constants provide a  measurement of  closeness of block
subspaces of $X$ to $\ell_1$. Clearly $X$ is asymptotic $\ell_1$ w.r.t.
$(e_i)$ if and only if $\delta_1(X)>0$.  The asymptotic $\ell_1$ constant
of $X$ is then equal to $\delta_1(X)^{-1}$. On the other hand we also
have
\begin{prop}\label{D.f}
  $X$ contains a subspace isomorphic to $\ell_1$ if and only if
  $\dot\delta_\alpha (X) >0$ for all $\alpha<\omega_1$.
\end{prop}
\begin{proof}
  This follows from Bourgain's $\ell_1$ index of a Banach space $X$ which
  we recall now.  For $0 < c < 1$, $\T(X,c)$ is the tree of all finite
  normalized sequences $(x_i)_1^k \subseteq X$ satisfying $\|\sum_1^k a_i
  x_i\|\ge c\sum_1^k |a_i|$ for $(a_i)_1^k\subseteq \real$.  The order on
  the tree is $(x_i)_1^k \le (y_i)_1^n$ if $k\le n$ and $x_i=y_i$ for
  $i\le k$.  For ordinals $\beta<\omega_1$ we define $\D^\beta (\T(X,c))$
  inductively by $\D^1 (\T(X,c)) = \{ (x_i)_1^k \in \T(X,c) :(x_i)_1^k$
  is not maximal$\}$.  $\D^{\beta+1} (\T(X,c)) = \D^1 (\D^\beta
  (\T(X,c)))$ and $\D^\beta (\T(X,c)) = \bigcap_{\gamma<\beta} \D^\gamma
  (\T(X,c))$ if $\beta$ is a limit ordinal.  The index $\I (X)$ is
  defined by $\I(X)=\sup_{0<c<1} \inf \{\beta : \D^\beta (\T(X,c))
  =\emptyset\}$, where the infimum is set equal to $\omega_1$ if no such
  $\beta$ exists.  Bourgain showed that for a separable space $X$,
  $\I(X)<\omega_1$ if and only if $X$ does not contain a subspace
  isomorphic to $\ell_1$ \cite{B}.
 
  Now observe that if $\F$ is a regular set of finite subsets of
  $\nat$ then $D(\F) = \{F\in \F: F\cup \{k\} \in \F$ for some
  $F<k\}$.  It follows that if $\delta_{\F}(x_i) >0$ for some basic
  sequence $(x_i)$ in $X$ then $\I(X)\ge I(\F)$.  Hence by
  Proposition~\ref{C.n}, if $\dot\delta_\alpha (X) >0$ for every
  $\alpha < \omega_1$ then $\I(X) = \omega_1$, hence $X$ contains a
  subspace isomorphic to $\ell_1$.  The converse implication is
  obvious.
\end{proof}

Other facts about Bourgain's $\ell_1$ index can be found in \cite{JO}.

The next lemma  collects some simple observations about the asymptotic
constants.
\begin{lemma}\label{D.e}
  Let $(e_i)$ be a basis for $X$ and let $(x_i)\prec (e_i)$.  Let $\F$
  and $\G$ be regular classes of finite subsets of $\nat$.
\begin{description}
\item[a)] $\delta_\F (e_i) \le   \delta_\F (x_i)$
and $\dot\delta_\F (x_i)\le \dot\delta_\F (e_i)$;
\item[b)] $\delta_\F (e_i) \le \dot\delta_\F (e_i) \le
  \ddot\delta_\F (e_i)$;
\item[c)] $\inf_n \delta_n (e_i) >0$ iff $({e_i}/{\|e_i\|})$ is
  equivalent to the unit vector basis of $\ell_1$;
\item[d)] $\ddot\delta_\F (e_i) = \sup_{(x_i)\prec (e_i)} \sup
  \{\dot\delta_\F ((x_i),|\cdot|) : |\cdot| \mbox{\rm \ is an equivalent
    norm on } [x_i]_{i\in\nat}\}$;
\item[e)] $\delta_{\F[\G]} (x_i)\ge \delta_\F (x_i) \delta_\G (x_i)$.
\end{description}
\end{lemma}
\begin{proof}
  a) and b) are immediate; the first part of a) uses that
  $\F(M)\subseteq \F$.  c) follows from the fact that
  $\bigcup_{n=1}^\infty \S_n$ contains all finite subsets of $\{2, 3,
  \ldots\}$.  d) is true because if $Y\subseteq X$ and $|\cdot|$ is an
  equivalent norm on $Y$ then $|\cdot|$ can be extended to an
  equivalent norm on $X$. For e) notice that if $(y_i)_i^k$ is
  $\F[\G]$-admissible w.r.t. $(x_i)$, then it can be blocked in a
  $\F$-admissible way into successive blocks each of which consists of
  $\G$-admissible vectors (w.r.t. $(x_i)$).  This directly implies the
  inequality.
\end{proof}

The most important situation for the study of the constants
$\delta_\alpha$ is when the whole sequence $(\delta_\alpha)_
{\alpha<{\omega_1}}$ is stabilized on a nested sequence of block
subspaces. This leads to the concept of the $\Delta$-spectrum of $X$
to be all possible stabilized limits of $\delta_\alpha$'s of block
bases. We formalize it in the following definition.

\begin{defin}\label{D.g}
  Let $X$ be a Banach space and let $\gamma=(\gamma_\alpha)_{\alpha
    <\omega_1}\subseteq \real$. We say that a basic sequence $(x_i)$
  in $X$ {\em \Dstabs\/} $\gamma$ if there exist
  $\ep_n\downarrow 0$ so that for every $\alpha <\omega_1$ there
  exists $m\in\nat$ so that for all $n\ge m$ if $(y_i) \prec
  (x_i)_n^\infty$ then $|\delta_\alpha (y_i)-\gamma_ \alpha| <\ep_n$.
 
  Let $X$ have a basis $(e_i)$.  The {\em $\Delta$-spectrum of $X$},
  $\Delta (X)$, is defined to be the set of all $\gamma$'s so that there
  exists $(x_i)\prec (e_i)$ such that $(x_i)$ \Dstabs\/
  $\gamma$.  By $\ddot\Delta (X)$ we denote the set of all $\gamma$'s so
  that $(x_i)$ \Dstabs\/ $\gamma$ for some $(x_i)\prec (e_i)$,
  under some equivalent norm $|\cdot|$ on $[x_i]_{i\in\nat}$.
\end{defin}

\begin{remark}\label{D.ga}\rm
  It is important to note that the asymptotic constants
  $\delta_\alpha(y_i)$ considered here and appearing in the definition
  of the spectrum $\Delta (X)$ refer to the admissibility with respect to
  the block basis $(y_i)$ itself. It is sometimes convenient, however,
  to consider asymptotic constants that keep a reference level for
  admissibility fixed when passing to block bases. Precisely, if
  $(e_i)$ is a basis in $X$ and $(x_i) \prec (e_i)$, we define
  $\delta_{\F} ((x_i),(e_i)) $ as the supremum of $\delta \ge 0$ such
  that whenever $(y_i)_1^k \prec (x_i)$ is $\F$-admissible w.r.t.
  $(e_i)$ then $\|\sum_1^k y_i\|\ge \delta \sum_1^k \|y_i\|$. Clearly,
  $\delta_\F (e_i) \le \delta_\F ((x_i),(e_i))\le \delta_\F (x_i)$. We
  can then define the spectrum $\Delta(X, (e_i))$ by replacing
  $\delta_\alpha (y_i)$ by $\delta_{\S_\alpha} ((y_i),(e_i))$, in
  Definition~\ref{D.g} above. Let us also note that it has been proved
  in \cite{AO} that these two concepts of spectrum actually coincide
  and $\Delta(X, (e_i))= \Delta(X)$.
\end{remark}

\begin{remark}\label{D.gb}\rm
  The definition of $\S_\alpha$  for $\alpha \ge \omega_0$ depended
  upon  certain choices made at limit ordinals. It follows that the
  constants $\delta_\alpha (e_i)$  also depend upon the particular
  choice  of $\S_\alpha $. However $\Delta (X)$  is independent
  of the choice of each $\S_\alpha$. Indeed, this follows from a
  consequence  of Propositions~\ref{C.n} and \ref{C.j}.
  If $\S_\alpha$  and $\bar\S_\alpha$  are two choices for the
  Schreier class then there exist subsequences of $\nat$, $M$ and  $N$
  such that $\S_\alpha (N) \subseteq \bar\S_\alpha$ and
  $\bar\S_\alpha (M) \subseteq \S_\alpha$. We also deduce that the
  constants $\dot \delta_\alpha$  and $\ddot \delta_\alpha$  are
  independent  of the particular choice of $\S_\alpha$.
\end{remark}

The following  stabilization argument shows that $\Delta (X)$ is always
non-empty.
\begin{prop}\label{D.h}
  Let $X$ be a Banach space with a basis $(e_i)$.  Then there exists
  $\gamma=(\gamma_\alpha)_{\alpha< \omega_1}$ and $(x_i) \prec (e_i)$
  so that $(x_i)$ \Dstabs\/ $\gamma$. In particular, $\Delta
  (X) \ne \emptyset$.
\end{prop}
\begin{proof}
  Fix $\ep_n \downarrow 0$.  If $[e_i]_{i\in\nat}$ contains $\ell_1$,
  then, since $\ell_1$ is not distortable, we can choose a normalized
  sequence $(x_i)\prec (e_i)$ with $\|\sum_n^\infty a_ix_i\| \ge
  (1-\ep_n) \sum |a_i|$ for all $(a_i)$; thus the proposition follows
  with $\gamma_\alpha =1$ for all $\alpha$.
 
  If $[e_i]_{i\in\nat}$ does not contain $\ell_1$ then by
  Proposition~\ref{D.f},
  $\dot\delta_\alpha (e_i)>0$ for at most countably many $\alpha$'s.
 
  Fix an arbitrary $\alpha < \omega_1$. It follows from
  Lemma~\ref{D.e} that if $(y_i) \prec (e_i)$ then $\dot\delta_\alpha
  ((y_i)_n^\infty) = \dot\delta_\alpha (y_i)$ for all $n$.  Since
  $\dot\delta_\alpha (y_i)\le \dot\delta_\alpha (z_i)$ whenever $(y_i)
  \prec (z_i)$,  by a standard  argument we can stabilize
  $\dot\delta_\alpha$.  That is, given $(w_i)\prec (e_i)$ we can find
  $(z_i)\prec (w_i)$ so that
  $$\gamma_\alpha \equiv \dot\delta_\alpha (z_i) = \dot\delta_\alpha
  (y_i) \quad\mbox{for all } (y_i)\prec (z_i)\ .$$ (To do this, construct
  $(w_i) \succ (z_i^{(1)}) \succ (z_i^{(2)})\succ \ldots$ such that
  $\dot\delta_\alpha (z_i^{(k+1)}) \le \inf\{\dot\delta_\alpha (y_i) :
  (y_i)\prec (z_i^{(k)})\} + 2^{-k}$, for every $k$, and set $z_i =
  z_i^{(i)}$ for all $i$.)
 
  Now choose by induction $(z_i) \succ (x_i^{(1)}) \succ (x_i^{(2)})\succ
  \ldots$ such that
  $$|\delta_\alpha(x_i^{(n+1)}) - \dot\delta_\alpha(x_i^{(n)})|=
  |\delta_\alpha(x_i^{(n+1)}) - \gamma_\alpha| \le \ep_n
  \quad\mbox{for all }n\ ,$$ and let $x_i = x_i^{(i)}$ for all $i$.
  Then $|\delta_\alpha ((x_i)_n^\infty) - \gamma_\alpha| <\ep_n$ for
  all $n$.  If $(y_i) \prec (x_i)_n^\infty$ then $\delta_\alpha
  ((x_i)_n^\infty) \le \delta_\alpha (y_i) \le \dot\delta_\alpha
  (y_i)= \gamma_\alpha$.
 
  Then using this and a diagonal argument for the countably many
  $\alpha$'s so that $\dot\delta_\alpha (e_i) >0$ we obtain the
  proposition.
\end{proof}

Our next proposition collects some  basic  facts about the $\Delta$-spectrum.

\begin{prop}\label{D.j}
  Let $X$ have a basis $(e_i)$.
\begin{description}
\item[a)]$\Delta (X) \ne \emptyset$ and if $\gamma \in\Delta (X)$ then
 $\gamma_\alpha \in [0,1]$ for $\alpha < \omega_1$.
\item[b)] $X$ contains $\ell_1$ iff there exists $\gamma\in\Delta (X)$
  with $\gamma_1=1$.
\item[c)] If $\gamma\in\Delta (X)$ then $\gamma_\alpha
  \ge\gamma_\beta$ if $\alpha \le\beta <\omega_1$.
\item[d)] If $\gamma\in\Delta (X)$ and $\alpha,\beta <\omega_1$, then
  $\gamma_\alpha\gamma_\beta \le \gamma_{\beta+\alpha}$.
\item[e)] If $\gamma\in\Delta(X)$ and $\alpha <\omega_1$, $n\in\nat$
  then $\gamma_{\alpha\cdot n} \ge (\gamma_\alpha)^n$.
\item[f)] If $\gamma\in \Delta (X)$ then $\gamma$ is a continuous
  function of $\alpha$.
\item[g)] $\ddot\delta_\alpha (X) = \sup \{\gamma_\alpha : \gamma\in
  \ddot\Delta (X)\}$.
\end{description}
\end{prop}

\begin{proof}
  We have already seen the  non-trivial part of a) and one implication in
  b).  Next, e) follows immediately from d) while f) and g) follow from
  the relevant definitions, using c) to get f).
 
  To complete b) note that if $\gamma_1=1$ then $\gamma_\alpha=1$, for
  all $\alpha < \omega_1$ (for $\alpha = \beta +1$ this follows from
  d) and for $\alpha$ a limit ordinal---from f)). Thus by
  Proposition~\ref{D.f}, $X$ contains $\ell_1$.
 
  c) Let $\gamma\in \Delta (X)$ and $\alpha \le\beta<\omega_1$.  For
  $n\in\nat$ let $\nat_n = (n,n+1,\ldots)$.  Let $(x_i)$ stabilize
  $\gamma$.  Given $m\in\nat$ choose $n \ge m $ by
  Proposition~\ref{C.b} so that $\S_\alpha (\nat_n) \subseteq \S_\beta
  (\nat_m)$.  It follows that $\delta_\alpha ((x_i)_n^\infty) \ge
  \delta_\beta ((x_i)_m^\infty)$.  Letting $m\to\infty$ we get
  $\gamma_\alpha \ge \gamma_\beta$.
 
  d) Let $(y_i)$ be basic.  By Proposition~\ref{C.b} there exists $M$
  with $\S_{\beta+\alpha} (M) \subseteq \S_\alpha [\S_\beta]$.  It
  follows that $\delta_{\beta+\alpha} ((y_i)_M) \ge \delta_{\S_\alpha
    [\S_\beta]} (y_i)$.  By Lemma~\ref{D.e} we see that
  $\delta_{\S_\alpha[\S_\beta]} (y_i) \ge \delta_\alpha (y_i)
  \delta_\beta (y_i)$.  Thus $\delta_{\beta+\alpha} ((y_i)_M) \ge
  \delta_\alpha (y_i)\delta_\beta (y_i)$.  Using this for $(y_i) =
  (x_i)_{i\in\nat_n}$ where $(x_i)$ stabilizes $\gamma$, we obtain
  that $\gamma_{\beta+\alpha} \ge\gamma_\alpha \gamma_\beta$.
\end{proof}
\begin{remark}\label{D.ja}\rm
  It is often useful to note that the constants $\delta_n$ satisfy
  conditions c) and d) for natural numbers. If $m, n \in \nat$ and $m \le
  n$ then $\delta_m (x_i) \ge \delta_n (x_i)$ and $\delta_{m+n} (x_i) \ge
  \delta_m (x_i)\delta_n (x_i)$, hence also $\delta_{mn} (x_i) \ge
  (\delta_n (x_i))^m$ (because $\S_m \subseteq \S_n$  and $\S_n[\S_m]
  = \S_{m+n}$).
\end{remark}

It is well known that the supermultiplicativity property d) of sequences
$\gamma \in \Delta (X)$ formally implies a ``sub-power-type'' behavior
of $\gamma$, which we shall find useful in various situations. This
depends on an elementary lemma. For two sequences $(b_n), (c_n)\subseteq
(0,1]$ we shall write $c_n \ll b_n$ to denote that $\lim_n b_n/c_n =
\infty$.
\begin{lemma}\label{D.l}
  Let $(b_n)\subseteq (0,1]$ satisfy $b_{n+m} \ge b_nb_m$ for all
  $n,m\in\nat$.  Then $\lim_n b_n^{1/n}$ exists and equals  $\sup_n
  b_n^{1/n}$. Moreover, for every $0 < \xi < \lim_n b_n^{1/n}$
  we have $\xi^n \ll b_n$.
\end{lemma}
\begin{proof}
  Let $a_n=\log (b_n^{-1})$.  Then $a_n\ge0$ and $a_{n+m} \le a_n +
  a_m$ for all $n,m$.  It suffices to prove that ${a_n}/n\to
  a\equiv \inf_m \{{a_m}/m\}$.  Given $\ep>0$ choose $k$ with
  $|{a_k}/k -a| <\ep$.  For $n>k$, ${a_n}/n - a< {a_n}/n -
  {a_k}/k +\ep$.  Setting $n=pk+r$, $0\le r<k$ and using $a_{p_k}
  \le pa_k$ we obtain
\begin{eqnarray*}
  \frac{a_n}n - \frac{a_k}k +\ep
  & \le & \frac{a_{pk} + a_r}n - \frac{a_k}k +\ep
  \le \frac{pa_k}{pk+r} + \frac{a_r}n - \frac{a_k}k +\ep\\
  & \le  & \frac{pa_k}{pk} - \frac{a_k}k + \frac{a_r}n +\ep =
  \frac{a_r}n +\ep\ .
\end{eqnarray*}
The first part of the lemma follows. The moreover part can be easily
proved by  contradiction.
\end{proof}

We have an  immediate  corollary.
\begin{cor}\label{D.la}
  Setting $\widehat\gamma_\alpha = \lim_n (\gamma_{\alpha\cdot n})^{1/n}$
  for $\alpha < \omega_1$ we have that for every $0 < \xi <
  \widehat\gamma_\alpha $, $\xi^n \ll \gamma_{\alpha\cdot n} \le
  \widehat\gamma_\alpha^{\,n}$, for all $\alpha < \omega_1$ and $n \in
  \nat$.
 
  Setting $\widehat\delta = \lim_n (\delta_{n}(x_i))^{1/n}$, for a basic
  sequence $(x_i)$, we have that for every $0 < \xi <\widehat\delta$,
  $\xi^n \ll \delta_n(x_i) \le \widehat\delta^n$ for all $n \in \nat$.
\end{cor}

There is an interesting connection between the constants
$\ddot\delta_\alpha(X)$ which allow for renormings of a given space
$X$, and the supermultiplicative behavior of $\gamma \in \Delta (X)$,
in particular of ${\widehat\gamma}_\alpha$, which involved the
original norm only.
\begin{prop}\label{D.m}
  Let $X$ have a basis $(e_i)$ and let $\gamma\in\Delta (X)$.  Then there
  exists $(y_i)\prec (e_i)$ so that $(y_i)$ \Dstabs\/ $\gamma$
  and so that for all $\alpha <\omega_1$, $\ddot\delta_\alpha (y_i) =
  \lim_n (\gamma_{\alpha\cdot n})^{1/n}\equiv \widehat\gamma_\alpha$.
\end{prop}

The argument is based on the following renorming result which we shall
use again.
\begin{prop}\label{D.n}
  Let $Y$ be a Banach space with a bimonotone basis $(y_i)$.  Let
  $\alpha<\omega_1$ and $n\in\nat$.  Then there exists an equivalent
  bimonotone norm $\norm{\cdot}$ on $Y$ with $\delta_\alpha
  ((y_i),\norm{\cdot}) \ge \Bigl(\delta_{[\S_\alpha]^n} (y_i)
  \Bigr)^{1/n}$.
\end{prop}
\begin{proof}
  Denote the original norm on $Y$ by $|\cdot|$  and set
  $\theta = \delta_{[\S_\alpha]^n} (y_i)$. For $0 \le j \le n$
  define  a norm $|\cdot|_j$ on $Y$ by
  \begin{eqnarray*}
    |y|_j = \sup \biggl\{ \theta^j \sum_1^\ell |E_iy| &:& (E_i y)_1^\ell
    \mbox{ is $[\S_\alpha]^j$-admissible w.r.t. } (y_i) \\  
    && \mbox{ and
    $E_1<\cdots<E_\ell$ are adjacent intervals} \biggr\}\ .
  \end{eqnarray*}
  Here we take $[\S_\alpha]^0 = \S_0$ so that $|y|_0 = |y|$.  For $0
  \le j \le n$ we have $|y|_j \ge \theta^j |y|$ and $|y| \ge
  \theta^{n-j} |y|_j$. The former inequality follows trivially from
  the definition of $|\cdot|_j$ and the latter from the
  fact that any $[\S_\alpha]^j$-admissible  family is
  $[\S_\alpha]^n$-admissible and the  definition of $\theta$.

  Set  $\norm{y} = \frac1n \sum_0^{n-1} |y|_j$  for $y \in
  Y$. Then $\norm{\cdot}$ is an equivalent norm on $Y$.
 
  Let $(x_s)_1^r$ be  $\alpha$-admissible w.r.t.  $(y_i)$.  First observe
  that $|\sum_1^r x_s| \ge \theta\sum_1^r |x_s|_{n-1}$.  Indeed,
  arbitrary $[\S_\alpha]^{n-1}$-admissible decompositions for each $x_s$
  can be put together to give a $[\S_\alpha]^{n}$-admissible decomposition
  for $\sum_1^r x_s$, thus the estimate follows from the definition of
  $|\cdot|_{n-1}$ and the fact that $\delta_{[\S_\alpha]^n} (y_i) =
  \theta^n$.  To be more precise,  for $1 \le s \le r$ choose adjacent
  intervals of integers $E_1^s<\cdots<E_{k(s)}^s$ so that
  $(E_j^s)_1^{k(s)}$ is $[\S_\alpha]^{n-1}$ admissible and
  $$|x_s|_{n-1} = \theta^{n-1}\sum_{j=1}^{k(s)} |E_j^s x_s|\ .$$
  Let $F_j^s = E_j^s$ if $j < k(s)$ and
  $F_{k(s)}^s = [\min E_{k(s)}^s, \min E_1^{s+1})$ if $s < r$ and
  $F_{k(r)}^r = E_{k(r)}^r$. Then
  $F_1<\cdots<F_{k(1)}^1<\cdots<F_{k(r)}^r$  are
  $[\S_\alpha]^n$-admissible  adjacent intervals of $\nat$ and so
  \begin{eqnarray*}
    \Bigl|\sum_{l=1}^r x_l \Bigr| &\ge& \theta^n
    \sum_{s=1}^r \sum_{j=1}^{k(s)} \bigl| F_j^s \bigl(\sum_{l=1}^r
    x_l\bigr)\bigr|\\
    &\ge& \theta \sum_{s=1}^r \theta^{n-1} \sum_{j=1}^{k(s)}|E_j^s (x_s)|
    = \theta \sum_{s=1}^r |x_s|_{n-1}\ ,
  \end{eqnarray*}
  (since $\bigl| F_{k(s)}^s \bigl(\sum_{l=1}^r x_l\bigr)\bigr| =
  | F_{k(s)}^s( x_s + x_{s+1})| \ge |E_{k(s)}^s (x_s)|$
  if $s < r$, using that the norm is monotone).

  Similarly, $|\sum_1^r x_s|_{j+1} \ge\theta\sum_1^r |x_s|_j$
  for $j=1,2,\ldots,n-2$, by the definitions of $|\cdot|_{j+1}$ and
  $|\cdot|_j$.  Thus
  $$ \Big|\!\Big|\!\Big| \sum_1^r x_s \Big|\!\Big|\!\Big| =\frac1n
  \sum_{j=0}^{n-1} \Big| \sum_{s=1}^r x_s\Big|_j \ge \frac{\theta}n
  \sum_{s=1}^r |x_s|_{n-1} + \frac1n \biggl( \sum_{j=0}^{n-2}\theta
  \sum_{s=1}^r |x_s|_j\biggr) \ .$$ Thus $\norm{\sum_1^r x_s}
  \ge\theta\sum_1^r \norm{x_s}$.
\end{proof}

\begin{prop}\label{D.na}
  Let $X$ be an asymptotic $\ell_1$ space and let $\gamma\in \Delta
  (X)$.  If $(e_i)\prec X$ \Dstabs\/ $\gamma$ then for all $\ep_i
  \downarrow 0$ there exists $(x_i)\prec (e_i)$ and an equivalent norm
  $|\cdot|$ on $[x_i]$ satisfying
  \begin{description}
  \item[a)] For all $n$ and $x \in \langle x_i\rangle_n^\infty$ we
    have
    $$ \|x|| \le |x| \le (2 + \ep_n)\|x\|\ .$$
  \item[b)] $(x_i)$ is bimonotone for $|\cdot|$.
  \item[c)] $(x_i)$ \Dstabs\/ $\bar\gamma \in \Delta(X,
    |\cdot|)$ with $\bar\gamma_\alpha \ge \gamma_\alpha$ for all
    $\alpha < \omega_1$.
  \end{description}
\end{prop}
\begin{proof}
  We may assume that $[e_i]$ does not contain $\ell_1$. Thus by
  Rosenthal's theorem \cite{R2} there exists $(x_i)\prec (e_i)$ which
  is normalized and weakly null. By passing to a subsequence of
  $(x_i)$ we may assume that for all $n < m$ and $(a_i)_1^m \subseteq
  \real$, $\|\sum_1^n a_i x_i \| \le (1 + \bar\ep_n) \|\sum_1^m a_i
  x_i \|$, where $\bar\ep_n = \ep_n/2$.

  Define the norm $|\cdot|$ for  $x \in X$  by
  $$|x| = \sup \{\|Ex\|: E \mbox{ is an interval}\}\ .$$
  Passing to a block basis of
  $(x_i)$ we may assume that $(x_i)$ \Dstabs\/ some
  $\bar\gamma \in \Delta(X, |\cdot|)$. For $x = \sum_{i=n}^m   a_i x_i$ 
  with  $|x| = \| Fx \|$  we have
  $$ \|x\| \le  |x| \le \|\sum_{i=n}^{\max F}   a_i x_i \|
      +  \|\sum_{i=n}^{\min F-1}   a_i x_i \|
      \le   2 (1 + \bar\ep_n)\|x\| = (2 +   \ep_n)\|x\|\ . $$
  Thus a) holds and b) is immediate. It remains to
  check c). Fix $\alpha < \omega_1$ and $m \in \nat$. Let $x_m <
  y_1<\cdots<y_\ell$ (w.r.t. $(x_i)_m^\infty$) where $(y_i)_1^\ell$ is
  $\alpha$-admissible w.r.t. $(x_i)_m^\infty$ and hence w.r.t.
  $(e_i)_m^\infty$. Choose intervals $E_1< \cdots < E_\ell$  such that
  $|y_i| = \|E_i y_i\|$ for $i \le \ell$  and $E_i \subseteq [\min
  \supp(y_i), \max \supp(y_i)]$. Define $(F_i)_1^\ell$ to be
  adjacent intervals so that $\min F_i = \min E_i$.
  Thus $F_i = [\min E_i, \min E_{i+1})\subseteq \nat$ for $i < \ell$
  and $F_\ell = E_\ell$. Let $F = \bigcup_1^\ell F_i$. Then, by
  Remark~\ref{D.ca},
  \begin{eqnarray*}
    \Bigl|\sum_1^\ell y_i\Bigr|&\ge& \Bigl\|F(\sum_1^\ell y_i)\Bigr\|
    \ge \delta_\alpha ((e_i)_m^\infty)\sum_{j=1}^\ell
    \|F_j(\sum_1^\ell y_i)\|\\
    &\ge& \delta_\alpha ((e_i)_m^\infty) (1 +\bar \ep_m)^{-1}
        \sum_{j=1}^\ell \|F_jy_j\| =
        \delta_\alpha ((e_i)_m^\infty)  (1 +\bar \ep_m)^{-1}\sum_{j=1}^\ell
        |y_j|\ .
  \end{eqnarray*}
  It follows that $\delta_\alpha((x_i)_m^\infty, |\cdot|) \ge
  \delta_\alpha ((e_i)_m^\infty)  (1 +\bar \ep_m)^{-1}$. Letting $m \to \infty$
  we obtain $\bar\gamma_\alpha \ge \gamma_\alpha$.
\end{proof}

\begin{remark}\label{D.nb}\rm
  It is worth noting the following.  Let $(e_i)$ be a basic sequence
  in $X$ $\Delta$-stabilizing $\gamma \in \Delta (X)$. Then there
  exists $(x_i) \prec (e_i)$ and an equivalent monotone norm $|\cdot|$
  on $[x_i]$ so that $(x_i)$ \Dstabs\/ $\gamma \in \Delta(X,
  |\cdot|)$. Furthermore $\Bigl||x| - \|x\|\Bigr| < \ep_n$ for
  $x\in\langle x_i\rangle_n^\infty$ and some $\ep_n \downarrow 0$.
  Assuming as we may that $[e_i]$ does not contain $\ell_1$, this is
  accomplished by taking $(x_i)$ to be a suitable weakly null block
  basis of $(e_i)$ and setting $|\sum a_i x_i| = \sup_n \|\sum_1^n a_i
  x_i\|$.
\end{remark}

A similar argument yields

\begin{prop}\label{D.nc}
  Let $\F$ be a regular set of finite subsets of $\nat$ and let
  $(e_i)$ be a basis for $X$. Given $\ep >0$ and $\ep_i \downarrow 0$
  there exists an equivalent norm $|\cdot|$ on some block subspace
  $[x_i]\subseteq X$ satisfying a) and  b) of Proposition~\ref{D.na}
  and $\delta_\F((x_i), |\cdot|) \ge \delta_\F (e_i) - \ep$.
\end{prop}

As a corollary to these propositions   we obtain
\begin{thm}\label{D.nd}
  Let $Y$ be a Banach space with a basis $(y_i)$. Let $\alpha <
  \omega_1$, $n \in \nat$,  $\ep >0$ and $\theta^n =
  \delta_{[\S_\alpha]^n}(y_i)$. Then there exists an equivalent norm
  $\norm{\cdot}$ on  $X = [x_i]\prec Y$ with $\delta_\alpha((x_i),
  \norm{\cdot}) \ge \theta - \ep$.
\end{thm}

\medskip \pf{Proof of  Proposition~\ref{D.m}}
  Let $(x_i) \prec (e_i)$ $\Delta$-stabilize $\gamma$ (for the original
  norm $\|\cdot\|$).
  We may assume that $X' = [x_i]$ does not contain $\ell_1$. 
  It follows that there
  exists $\alpha_0<\omega_1$ so that $\ddot\delta_\beta (X') =
  0=\gamma_\beta$ for all $\beta>\alpha_0$.
  Also  from  Lemma~\ref{D.e},   $\ddot\delta_\alpha (z_i) \le
  \ddot\delta_\alpha (w_i)$ if $(z_i)\prec (w_i)\prec (e_i)$; moreover,
  $\ddot\delta_\alpha ((z_i)_n^\infty)= \ddot\delta_\alpha (z_i)$ for all
  $n \in \nat$. We can therefore stabilize the $\ddot\delta_\alpha$'s (as
  in the proof of Proposition~\ref{D.h}) to find $(y_i)\prec (x_i)$ so
  that for all $\alpha \le\alpha_0$,  $\ddot\delta_\alpha (y_i) =
  \ddot\delta_\alpha (z_i)$ if $(z_i) \prec (y_i)$.  Of course $(y_i)$
  still \Dstabs\/ $\gamma$.  We shall prove that
  $\ddot\delta_\alpha (y_i) = \lim_n (\gamma_{\alpha\cdot n})^{1/n}$.
 
  Note that if $|\cdot|$ is an equivalent norm on $[y_i]$ and
  $\bar\gamma \in \Delta ((y_i),|\cdot|)$ then $\lim_n (\bar \gamma_
  {\alpha\cdot n})^{1/n} = \lim_n (\gamma_{\alpha\cdot n})^{1/n}$.
  Indeed if $(z_i) \prec (y_i)$ \Dstabs\/ $\bar\gamma$ in $|\cdot|$
  then since $(z_i)$ \Dstabs\/ $\gamma$ in $\|\cdot\|$ and the norms
  are equivalent, we obtain $c\bar\gamma_\beta \le \gamma_\beta \le
  d\bar\gamma_\beta$ for all $\beta<\omega_1$ and for some constants
  $c$, $d>0$.  Thus
  $$\bar\gamma_\alpha \le \sup_n (\bar\gamma_{\alpha\cdot n})^{1/n} =
  \lim_n (\gamma_{\alpha\cdot n})^{1/n}\ .$$ By Proposition~\ref{D.j}
  we obtain that $\ddot\delta_\alpha (y_i) \le \lim_n
  (\gamma_{\alpha\cdot n})^{1/n}$.
 
  Fix $\theta < \lim_n (\gamma_{\alpha\cdot n})^{1/n}$.  Thus there
  exists $n_0$ with $\theta^{n_0} <\gamma_{\alpha\cdot n_0}$.  Choose
  $(z_i) \prec (y_i)$ with $\theta^{n_0} <\delta_{\alpha\cdot n_0}
  (z_i)$.  By Corollary~\ref{C.d} there exists $M$ so that
  $[\S_\alpha]^{n_0}(M)\subseteq \S_{\alpha\cdot {n_0}}$, which yields
  $\delta_{\alpha\cdot n_0}(z_i) \le \delta_{[\S_\alpha]^{n_0}}
  ((z_i)_M)$.  So letting $(w_i) = (z_i)_M$ we have
  $\delta_{[\S_\alpha]^{n_0}} (w_i) > \theta^{n_0}$.  By
  Theorem~\ref{D.nd} there exists an equivalent norm $\norm{\cdot}$ on
  $[w_i']_{\nat}$, for some $(w_i')\prec (w_i)$ with $\delta_\alpha
  ((w_i'),\norm{\cdot}) >\theta$.  The reverse inequality,
  $\ddot\delta_\alpha (y_i) \ge \lim_n (\gamma_{\alpha\cdot
    n})^{1/n}$, follows. 
\endpf

As we will see in later sections, some further regularity properties
of sequences $\gamma \in \Delta (X)$ are closely related to distortion
properties of the space $X$, and they may or may not hold in general.
In contrast, the sequences $(\ddot\delta_\alpha)$ which allow for
renorming display a complete power type behavior. In fact, we will
give a comprehensive  description of behavior of such sequences in
Theorem~\ref{F.c} below.

In the result that follows we shall be particularly interested in part
c).
\begin{prop}\label{D.k}
  Let $X$ have a basis $(e_i)$.  Let $\alpha<\omega_1$ and $n\in\nat$.
\begin{description}
\item[a)] $\ddot\delta_{[\S_\alpha]^n} (X) = (\ddot\delta_\alpha(X))^n$
\item[b)] $\ddot\delta_{[\S_\alpha]^n}(X) = \ddot\delta_{\alpha\cdot
    n}(X)$
\item[c)] $\ddot\delta_{\alpha\cdot n} (X) = (\ddot\delta_\alpha
  (X))^n$
\end{description}
\end{prop}
\begin{proof}
  c) will follow from a) and b).
 
  a) Since for any equivalent norm $|\cdot|$ on $X$ we have
  $\delta_{[\S_\alpha]^n} ((y_i),|\cdot|) \ge {(\delta_\alpha
    ((y_i),|\cdot|))^n}$ (Lemma~\ref{D.e}, e)), the inequality
  $\ddot\delta_{[\S_\alpha]^n}(X) \ge (\ddot\delta_\alpha (X))^n$
  follows from g) of proposition~\ref{D.j}.  To see the reverse
  inequality let $|\cdot|$ be an equivalent norm on $X$, and let
  $(y_i)\prec (e_i)$ and $\theta >0$ satisfy $\delta_{[\S_\alpha]^n}
  ((y_i),|\cdot|) > \theta^n$.  By Theorem~\ref{D.nd} there exist
  $(x_i) \prec (y_i)$ and an equivalent norm $\norm{\cdot}$ on
  $[x_i]_{i\in\nat}$ such that $\delta_\alpha ((x_i),\norm{\cdot})>
  \theta$.  This completes the proof.
 
  b) As we have shown earlier, whenever $(y_i)\prec (e_i)$ and
  $|\cdot|$ is an equivalent norm, by Corollary~\ref{C.d}  there exists
  a subsequence $M$ such that $\delta_{[\S_\alpha]^n}
  ((y_i)_M,|\cdot|) \ge \delta_{\alpha\cdot n} ((y_i),|\cdot|)$.  It
  follows that $\ddot\delta_{[\S_\alpha]^n} (X) \ge
  \ddot\delta_{\alpha\cdot n} (X)$.  The reverse inequality follows by
  choosing $N$ with $\S_{\alpha\cdot n}(N) \subseteq [\S_\alpha]^n$.
\end{proof}

Let us introduce the following natural and convenient definition.
\begin{defin}\label{F.a}
  Let $X$ be an asymptotic $\ell_1$ space.  The {\em spectral index
    of\/} $X$, $I_\Delta (X)$, is defined to be
  $$I_\Delta (X) = \inf \{\alpha <\omega : \ddot\delta_\alpha (X) <1\}\ .$$
\end{defin}

\begin{thm}\label{F.c}
  If $X$ is an asymptotic $\ell_1$ space not containing $\ell_1$, then
  $I_\Delta(X) = \omega^\alpha$ for some $\alpha<\omega_1$.  If
  $I_\Delta(X) =\alpha_0$ and $\ddot\delta_{\alpha_0} (X)=\theta$ then
  $\ddot\delta_{\alpha_0\cdot n+\beta} (X) = \theta^n$ for all $n\in\nat$
  and $\beta <\alpha_0$. Finally, $\ddot\delta_{\beta} (X) = 0 $ for all
  $\alpha_0\cdot\omega \le \beta < \omega_1$.
\end{thm}
\begin{proof}
  For the proof of the first statement, it suffices to show that if
  $\beta <I_\Delta(X)$ then for all $n\in\nat$, $\beta\cdot n <I_\Delta
  (X)$ (\cite{M}, Thm.~15.5).  But by Proposition~\ref{D.k},
  $\ddot\delta_ {\beta\cdot n} (X) = (\ddot\delta_\beta (X))^n = 1$, so
  $\beta\cdot n <I_\Delta (X)$.
 
  Now let $\alpha_0= \omega^{\alpha}$ for some $\alpha$ and assume
  that $\ddot\delta_{\alpha_0} (X) = \theta$ for some $0 < \theta <1$.
  Fix $\beta < \alpha_0$. We first show that for any $ \ep > 0 $ we
  can find $(y_i)\prec X$ and an equivalent norm $\norm{\cdot}$ on
  $[y_i]_{\nat}$ with $\delta_\beta ((y_i),\norm{\cdot}) > 1-\ep$ and
  $\delta_{\alpha_0} ((y_i),\norm{\cdot}) > \theta-\ep$.  Indeed, let
  $\theta' = \theta - \ep$ and choose by Proposition~\ref{D.na}
  $(x_i)\prec X$ and an equivalent bimonotone norm $|\cdot|$ on $X$ so
  that $\delta_{\alpha_0} ((x_i),|\cdot|) >\theta'$.  Given $m\in\nat$
  we can choose a subsequence $N$ of $\nat$ so that $\S_{\alpha_0}
  \Bigl[[\S_\beta]^j \Bigr](N) \subseteq \S_{\beta\cdot j + \alpha_0} =
  \S_{\alpha_0}$ for $j=0,1,\ldots,m$; this follows from
  Proposition~\ref{C.b}, Corollary~\ref{C.d} and the fact that
  $\beta\cdot m + \omega^{\alpha} = \omega^{\alpha}$.  Let
  $(y_i)=(x_i)_N$ and $a^m\equiv \delta_{[\S_\beta]^ m}
  ((y_i),|\cdot|)$.  Note that since  $[\S_\beta]^m (N) \subseteq
  \S_{\alpha_0}$
  then $a^m\ge \theta'$ and so $a\ge (\theta')^{1/m}$.  For $y\in [
  y_i]_\nat$ and $0\le j\le m$ set
  \begin{eqnarray*}
  |y|_j = \sup \{a^j \sum_1^\ell |E_iy| &:& (E_iy)_1^\ell \mbox{ is
    $[\S_\beta]^ j$-admissible w.r.t. } (y_i) \\
    && \mbox{ and     $E_1<\cdots<E_k$  are adjacent intervals} \}\ .
  \end{eqnarray*}
 
  It can be checked by a straightforward calculation, using the choice
  of $N$ and that $(y_i)$ is monotone for $|\cdot|$, that
  $\delta_{\alpha_0} ((x_i),|\cdot|_j) \ge \delta_{\alpha_0}
  ((y_i),|\cdot|) >\theta'$ for $j=0,\ldots,m$.  For $y\in [
  y_i]_\nat$ set $\norm{y} = \frac1m \sum_{j=0}^{m-1} |y|_j$.  Then
  $\delta_{\alpha_0} ((y_i),\norm{\cdot}) >\theta'$ and from the proof
  of Proposition~\ref{D.n}, $\delta_\beta ((y_i),\norm{\cdot}) \ge
  a>(\theta')^{1/m}$. Taking $m$ such that $(\theta')^{1/m} \ge 1 -
  \ep$ we get what we wanted.
 
  Now by Proposition~\ref{C.b} there exists a subsequence $M$ of $\nat$
  with $\S_{\alpha_0+\beta} (M) \subseteq \S_\beta [\S_{\alpha_0}]$.  It
  follows that
  $$\delta_{\alpha_0+\beta} ((y_i)_M,\norm{\cdot}) > (1-\ep)\theta' =
  (1-\ep)(\theta-\ep)\   .$$
  Hence $\ddot\delta_{\alpha_0+\beta} (X) = \theta$. 
 
  The case of general $n$ is proved similarly, replacing $\alpha_0$ by
  $\alpha_0\cdot n$ above and recalling (Proposition~\ref{D.k}) that
  $\ddot\delta_{\alpha_0\cdot n} (X) = (\ddot\delta_{\alpha_0} (X))^n$.
  The last statement is obvious.
\end{proof}

\section{Examples--Tsirelson Spaces}

Our primary source of examples of asymptotic $\ell_1$ spaces with various
behaviors of asymptotic constants is the class of  mixed Tsirelson spaces
introduced by Argyros and Deliyanni in \cite{AD}.
\begin{defin}\label{D.b}
  Let $I\subseteq \nat$ and for $n\in I$ let $\F_n$ be a regular
  family of finite subsets of $\nat$.  Let $(\theta_n)_{n\in I}
  \subseteq (0,1)$ satisfy $\sup_{n\in I} \theta_n<1$.  The {\em mixed
    Tsirelson space\/} $T(\F_n,\theta_n)_{n\in I}$ is the completion
  of $c_{00}$ under the implicit norm
  $$\|x\| = \max \Biggl(\|x\|_\infty,  \sup_{n\in I} \sup \biggl\{ \theta_n
  \sum_{i=1}^k \|E_i x\| :(E_i)_{i=1}^k \mbox{ is
    $\F_n$-admissible}\biggr\}\Biggr)\ .$$
\end{defin}
It is shown in \cite{AD} that such a norm exists.  It is also proved
that if $I$ is finite or if $\theta_n \to 0$, then
 $T(\F_n,\theta_n)_{n\in I}$ is a reflexive Banach space, in which
the standard unit vectors $(e_i)$ form a 1-unconditional basis.  In
\cite{AD} it is proved that for an appropriate choice of $\theta_n$
and $\F_n$ the space $T(\F_n,\theta_n)_{n\in \nat}$ is arbitrarily
distortable. Deliyanni and Kutzarova \cite{DK} proved a result that
illustrates the possible complexity these spaces can possess. They
proved that a mixed Tsirelson space may uniformly contain
$\ell_\infty^n$'s in all subspaces.   Notice  that
the Tsirelson space $T$ satisfies $T = T(\S_n,2^{-n})_{n\in\nat}
=T(\S_1,2^{-1})$. For $0 < \theta < 1$ we denote the
$\theta$-Tsirelson space by $T_\theta = T(\S_1,\theta)$.

\begin{thm}\label{E.a}
  Let $(e_i)$  denote  the unit vector basis for $T$.
\begin{description}
\item[a)] If $(x_i)\prec (e_i)$ then for all $n$, $\delta_n(x_i)
  =2^{-n}$  and $\ddot\delta_n(x_i) =2^{-n}$.
\item[b)] For all $\gamma \in \Delta (T)$, $\gamma_n=2^{-n}$ for $n\in
  \nat$ and $\gamma_\alpha=0$ for $\alpha\ge \omega$.
\item[c)] For all $\gamma\in\ddot\Delta (T)$, $\gamma_n\le 2^{-n}$ for
  $n\in\nat$.
\item[d)] $I_\Delta (X) =1$  for all $X \prec T$.
\end{description}
\end{thm}
\begin{remark}\label{E.ac}\rm
  Condition a) immediately implies that for an arbitrary equivalent norm
  $|\cdot|$ on $T$ and  $(x_i) \prec (e_i)$, we have
  $\delta_1((x_i), |\cdot|)\le 1/2$. Since the asymptotic $\ell_1$
  constant is equal to $\delta_1^{-1}$, this improves the  constant in
  Proposition~\ref{B.e} from $\sqrt 2$ to 2.
\end{remark}
\begin{remark}\label{E.aa}\rm
  For $T_\theta$ we have $\delta_n (T_\theta)= \ddot \delta_n (T_\theta)
  = \theta^n$ for $n \in \nat$; and all other equalities and inequalities
  from Theorem~\ref{E.a} hold with appropriate modifications. Also,
  clearly, $I_\Delta (T_\theta)=1$.
\end{remark}
\medskip \pf{Proof of  Theorem~\ref{E.a}}
  a) By definition of the norm $\|\cdot\|$ for $T$, $\delta_n(e_i) \ge
  2^{-n}$ and so if $(x_i)\prec (e_i)$ then $\delta_n(x_i)\ge 2^{-n}$ as
  well.
 
  We next show that there exists $C<\infty$ so that $\delta_m(x_i) \le
  C2^{-m}$ for all $m$.  This will yield the equality for $\delta_n$.
  Indeed if for some $n$, $\delta_n(x_i)= A/2^n$ where $A>1$ then since
  $\delta_{nk} (x_i) \ge (\delta_n(x_i))^k$ (Remark~\ref{D.ja}),
  we would have that $C2^{-nk}\ge \delta_{nk}(x_i)\ge A^k2^{-nk}$ for all
  $k$, which is impossible.
 
  First we consider the case $(x_i)=(e_i)_{i\in M}$ where $M$ is a
  subsequence of $\nat$.  Let $\ep>0$, $n\in \nat$ and let $x= \sum_{i
    \in F} a_i e_i$ be an $(n,n-1,\ep)$-average of $(e_i)_{i\in M}$ (see
  Proposition~\ref{C.l} and Notation~\ref{C.m}).  Thus $\|x\|\ge 2^{-n}$.
  Iterating the definition of the norm in $T$ yields that $\|x\| =
  \sum_{i=1}^n 2^{-i} \sum_{j\in F_i} a_j$ where $(F_i)_{i=1}^n$
  partitions $F$ into sets with $F_i\in \S_i$ for $i\le n$.  Thus if $\ep
  <2^{-n}$,
  $$\|x\| = \|\sum_{i\in F} a_i e_i\| \le \sum_{i=1}^{n-1} 2^{-i}\ep
  + 2^{-n} \sum_{j\in F_n} a_j \le 2/2^n
  =  2/2^n \sum_{i\in F} \|a_i e_i\|  \ .$$
  Hence $\delta_n((e_i)_M)\le 2/2^n$.
 
  If $(x_i)$ is normalized with $(x_i) \prec (e_i)$ then by \cite{CJT}
  (see also \cite{CS}), there exists a subsequence $M$ such that $(x_i)$
  is $D$-equivalent to $(e_i)_{i\in M}$, where $D$ is an absolute
  constant (we let $m_i=\min \supp (x_i)$, and then $M=(m_i)$).  Thus
  $ \delta_n (x_i) \le D \delta_n ((e_i)_M) \le  2D/2^n $.
 
  To get the equality for $\ddot\delta_n$ we first observe that for any
  equivalent norm $|\cdot|$ on $T$ there is a constant $C'$ (depending on
  $|\cdot|$) such that $ \delta_n ((x_i),|\cdot|) \le C' \delta_n (x_i)$,
  and then we follow the previous argument.
 
  b) is immediate from the first part of a); and c) and d) follow from
  the second part of a). 
\endpf

\begin{remark}\label{E.ab}\rm
  For the subsequence $M= (m_i)$ above one could take any $m_i \in \supp
  (x_i)$ for all $i$.  In the space $T_\theta$, any normalized block
  basis is $D$-equivalent to $(e_i)_M$ as well, with the equivalence
  constant $D = c\theta^{-1}$, where $c$ is an absolute constant. The
  choice of a subsequence $M$ is the same as indicated above (for
  $\theta = 1/2$).
\end{remark}

The next example illustrates Theorem~\ref{F.c}.
\begin{example}\label{F.d}
  Let $\alpha<\omega_1$ and let $X= T(\S_{\omega^\alpha},\theta)$.
  Then
\begin{description}
\item[a)] $\ddot\delta_{\omega^\alpha\cdot n} (X) = \theta^n$ for
  $n\in\nat$
\item[b)] $I_\Delta (X) = \omega^\alpha$
\end{description}
\end{example}
\begin{proof}
  a) Let $(x_i)\prec X$ be a normalized block basis that
  \Dstabs\/  $\gamma\in\Delta (X)$.  Let $n\in\nat$ and
  let $\ep >0$.  Choose $N$ by Corollary~\ref{C.d} so that
  $[\S_{\omega^\alpha}]^n \supseteq \S_{\omega^\alpha\cdot n}(N)$
  and also $[\S_{\omega^\alpha}]^{n-1}(N)\subseteq
  \S_{\omega^\alpha\cdot(n-1)}$.
  Choose $x= \sum_F a_i x_i$ to be an $(\omega^\alpha\cdot
  n,\omega^\alpha \cdot (n-1),\ep)$ average of $(x_i)_N$ w.r.t.
  $(e_i)$, the unit vector basis of $X$.  Clearly $\|x\|\ge \theta^n$.
  As in $T$, $\|x\|$ is calculated by a tree of sets where the first
  level of sets is $\S_{\omega^\alpha}$-admissible, the second level is
  $[\S_{\omega^\alpha}]^2$-admissible and so on.
 
  If we stop this tree after $n-1$ levels, discarding sets which
  stopped before then and shrinking those sets which split the support
  of some $x_i$ we obtain for some $(E_ix)_1^\ell$ being
  $\omega^\alpha\cdot (n-1)$-admissible,
  $$\|x\| \le \theta^{n-1} \sum_1^\ell \|E_i x\| + \ep\ .$$ The next
  level of splitting may indeed split the supports of some of the
  $x_j$'s.  However since those $x_j$'s have not yet been split the
  contribution of $a_j x_j$ to the next level of sets is at most $a_j
  \theta^{-1}$.  Thus we obtain
  $$\|x\| \le \theta^n \Big( \sum a_j \theta^{-1} \Big) + \ep =
  \theta^{n-1} + \ep\ .$$ It follows that $\gamma_{\omega^\alpha\cdot
    n} \le \theta^{n-1} = \frac1{\theta}(\theta^n)$.
 
  Thus, just as in the case of $T$, $\gamma_{\omega^\alpha\cdot n} =
  \theta^n$.  Indeed, if $\gamma_{\omega^\alpha\cdot n_0} >
  \theta^{n_0}$ then
  $$\gamma_{\omega^\alpha\cdot n_0k} \ge (\gamma_{\omega^\alpha\cdot
    n_0} )^k > \frac1{\theta} \theta^{n_0k}$$ for large enough $k$
  (Proposition~\ref{D.j}), which is a contradiction.
 
  Similarly if $\gamma\in\ddot\Delta (X)$ then for some $C$,
  $\gamma_{\omega^\alpha\cdot n} \le C\theta^n$ and so
  $\gamma_{\omega^\alpha\cdot n} \le \theta^n$ for all $n$.  This
  yields that $\ddot\delta_{\omega^\alpha\cdot n} (X) = \theta^n$.
 
  b) The argument in Proposition~\ref{F.c}(b) yields this result: for
  $\beta <\omega^\alpha$ and $\ep >0$ there exists $(x_i)$ and
  $\norm{\cdot}$ with $\delta_\beta ((x_i),\norm{\cdot}) >1-\ep$.
\end{proof}

Before we pass to further examples, let us note a fundamental and useful
connection between the spectrum $\Delta(X)$ and a lower estimate for the
norm on some block subspace.

\begin{prop}\label{E.c}
  Let $X$ be an asymptotic $\ell_1$ space and let $(z_i)\prec X$ be a
  normalized bimonotone block basis $\Delta$-stabilizing some
  $\gamma\in\Delta(X)$ with $0<\gamma_1<1$.  Let $(e_i)$ be the unit
  vector basis of $T_{\gamma_1} \equiv T(\S_1,\gamma_1)$.  Then for
  all $\ep>0$ there exists a subsequence $(x_i)$ of $ (z_i)$
  satisfying for all $(a_i)\subseteq \real$
  $$\|\sum a_ix_i\| \ge (1-\ep)\|\sum a_ie_i\|_{T_{\gamma_1}}\ .$$
\end{prop}
\begin{proof}
  We shall prove the proposition in the case where $\gamma_1 = 1/2$
  (and  so $T_{\gamma_1} = T$).
  We shall describe below the argument in a general case, but the reader
  is advised to first test the special case when $\delta_1 (z_i) = 1/2$
  (when $\ep_n =0$ for all $n$ and the $m_i$'s can be omitted.)
  Choose integers $m_i\uparrow \infty$ so that $\sum_1^\infty
  2^{-{m_i}}< \ep$ and
  then choose $\ep_n\downarrow 0$ to satisfy, for all $k\in\nat$,
  \begin{eqnarray} \label{E.eqA}\nonumber   
  \prod_1^k \Bigl(\frac12 -\ep_{n(i)}\Bigr) & > &(1-\ep)2^{-k} \mbox{
    whenever } (n(i))_{i=1}^k \subseteq \nat \\
  &&\mbox{ satisfy for every $j$, } |\{i:n(i)=j\}| \le m_j\ .
  \end{eqnarray}
  Let $(x_i)$ be a subsequence of $(z_i)$ which satisfies: for all $n$,
  if $x_n\le y_1<\cdots <y_n$ w.r.t. $(x_i)$ then $\|\sum_1^n y_i\| >
  (\frac12 -\ep_n) \sum_1^n \|y_i\|$.  Such a sequence exists since
  $(z_i)$ \Dstabs\/ $\gamma$ with $\gamma_1 =1/2$.
 
  Let $x=\sum_1^\ell a_ix_i$ and assume that $\|\sum_1^\ell a_i e_i\|_T
  =1$.  We shall show that $\|x\| >(1-\ep)^2$.  If $\|\sum a_ie_i\|_T =
  |a_j|$ for some $j$ then $\|x\| = 1 $.  Otherwise for some
  1-admissible family of sets, $\|\sum a_ie_i\|_T = \frac12 \sum_{j=1}^n
  \|E_j (\sum a_ie_i)\|_T$.  Accordingly we have that (here is where
  the bimonotone assumption is used)
  $$\|x\| > \left(\frac12 -\ep_i\right) \sum_{j=1}^n \|E_j x\|$$
  where $i=\min (\supp E_1x)$.  We then repeat the step above for each
  $E_j x$.  Ultimately we obtain for some $J\subseteq \nat$,
  $$1=\|\sum a_ie_i\|_T = \sum_{i\in J} 2^{-\ell (i)} |a_i|$$ where
  $\ell(i)=$ the number of splittings before we stop at $|a_i|$.  We
  follow the same tree of splittings in getting a lower estimate for
  $\|x\|$ with one additional proviso.  Each splitting of $Ex$ in
  $\langle x_i\rangle$ will introduce a factor of $(\frac12-\ep_n)$ for
  some $n$.  A given factor $(\frac12 -\ep_n)$ may be repeated a number
  of times.  If any $(\frac12-\ep_n)$ is repeated ${m_n}$ times we shall
  discard the corresponding set $\|Ex\|$ at that instant.  By virtue of
  (\ref{E.eqA}) we thus obtain that $\|x\| \ge \sum_{i\in I}
  (1-\ep)2^{-\ell (i)}|a_i|$ where $I\subseteq J$ and $a_ix_i$ belonged
  to a discarded set for $i\in J\setminus I$.  However the contribution
  of the discarded sets to $\|\sum a_i e_i\|_T$ is at most
  $\sum_{n=1}^\infty 2^{-{m_n}} <\ep$ since from our construction for any
  given $n$ (where $(\frac12-\ep_n)$ is repeated ${m_n}$ times) we will
  discard at most one set, something of the form $2^{-k}\|Ex\|_T$ where
  $k\ge {m_n}$.  It follows that $\|x\| > (1-\ep)(\|x\|_T -\ep) =(1-\ep)^2$.
\end{proof}

The proof also yields the following block result.
\begin{cor}\label{E.g}
  Let $(z_i)$ be a bimonotone basic sequence in a Banach space $X$
  which \Dstabs\/ $\gamma \in \Delta (X)$ where $0<\gamma_1<1$.  Let
  $(e_i)$ be the unit vector basis of $T_{\gamma_1}$.  Then for all
  $\ep>0$ there exists a subsequence $(x_i)$ of $ (z_i)$ satisfying
  for all $(y_j)_1^k \prec (x_i)$ if $m_j = \min (\supp (y_i))$ w.r.t.
  $(x_i)$ then
  $$\Bigl\|\sum_1^k y_i\Bigr\| \ge (1-\ep) \Big\| \sum_1^k  \|y_j\|
       e_{m_j}\Big\|_{T_\gamma}\ .$$
\end{cor}

\begin{remark}\label{E.ga}\rm
  We can remove the bimonotone assumption on the norm if we have that
  for some $\ep_n \downarrow 0$, $\|y_0 + \sum_1^m y_i\| \ge (\gamma_1
  - \ep_n) \sum_1^m \|y_i\|$, whenever $z_n \le y_0 \le z_m < y_1 <
  \cdots < y_m$. Without either this assumption or the bimonotone
  property we obtain a slightly weaker result.
\end{remark}

\begin{thm}\label{E.gb}
  Let $X$ be an asymptotic $\ell_1$ space and let $(z_i) \prec X$ be a
  basic sequence $\Delta$-stabilizing some $\gamma \in \Delta(X)$,
  with $0 < \gamma_1 <1$. Then for all $\ep >0$ there exists a
  normalized $(x_i) \prec (z_i)$ satisfying for all $(a_i) \subseteq
  \real$
  $$\|\sum a_i x_i\| \ge \frac12 (1 - \ep) \|\sum a_i
  e_i\|_{T_{\gamma_1}}\ .$$
  Moreover if $(y_i)_1^k \prec (x_i)$ with $m_j = \min(\supp (y_i))$
  w.r.t. $(x_i)$ then one has
  $$\|\sum_1^k y_i\|\ge \frac12 (1 - \ep)
  \Bigl\|\sum_1^k \|y_i\| e_{m_i} \Bigr\|_{T_{\gamma_1}}\ .$$
\end{thm}
\begin{proof}
  By Proposition~\ref{D.na} there exists a $\|\cdot\|$-normalized
  $(x_i)\prec(z_i)$ and a bimonotone norm $|\cdot|$ on $[x_i]$ with
  $\|x\| \le |x|\le (2 +\ep) \|x\|$ for $x \in [x_i]$ and such that
  $(x_i)$ \Dstabs\/ $\bar\gamma \in \Delta (X, |\cdot|)$ with
  $\bar\gamma_1 \ge \gamma_1$. We may thus assume that $(x_i)$
  satisfies the conclusion of Corollary~\ref{E.g} for $|\cdot|$ and
  $\ep' $ such that $(1-\ep')/(2+\ep')=\frac12 (1-\ep)$.  Thus if
  $(y_i)_1^k$ is as in the statement of the theorem,
$$  \Bigl\|\sum_1^k y_i\Bigr\|\ge \frac{1}{2+\ep'}  \Bigl|\sum_1^k
  y_i\Bigr| \ge  \frac{1-\ep'}{2+\ep'}\Bigl\|\sum_1^k
   |y_i| e_{m_i} \Bigr\|_{T_{\bar\gamma_1}}
     \ge \frac12 (1-\ep) \Bigl\|\sum_1^k
   \|y_i\| e_{m_i} \Bigr\|_{T_{\gamma_1}}\ .$$
\end{proof}

\medskip
The following can be proved by an argument similar to that in
Proposition~\ref{E.c}.
\begin{prop}\label{E.gc}
  Let $X$ be an asymptotic $\ell_1$ space and let $(z_i)\prec X$ be a
  normalized bimonotone block basis $\Delta$-stabilizing $\gamma \in
  \Delta (X)$. Let $\alpha < \omega_1$ with $0 < \gamma_\alpha <1$ and
  let $\ep >0$. Then there exists a subsequence $(x_i)$ of $(z_i)$
  satisfying the following: if $(y_i)_1^k \prec (z_i)$ with
  $\min(\supp(y_i))=m_i $ (w.r.t. $(x_i)$) then
  $$ \Bigl\|\sum_1^k y_i\Bigr\| \ge (1-\ep)  \Bigl\|\sum_1^k
   \|y_i\| e_{m_i} \Bigr\|_{T({\S_\alpha}, {\gamma_\alpha})}\ .$$
\end{prop}

The next example is a space $X$ for which the sequences of asymptotic
constants $(\delta_\alpha(X))$ and $(\ddot\delta_\alpha(X))$ are
``essentially'' the same as for Tsirelson's space $T$; still, $X$ and
$T$ have no common subspaces--no subspace of $X$ is isomorphic to a
subspace of $T$.  It is worth noting that $X$ also has the  property that
the sequence $\ddot\delta = (\ddot\delta_\alpha(X))$ does not belong
to $\ddot\Delta(X)$.

\begin{example}\label{E.h}
  Let $0<c<1$ and let $X= T(\S_n ,{c}{2^{-n}})_{n\in\nat}$.  Then
\begin{description}
\item[a)] $\ddot\delta_n(X) = 2^{-n}$ for all $n$
\item[b)] For all $\gamma\in \ddot\Delta (X)$, $\gamma_n <2^{-n}$ for
  all $n$.
\item[c)] No subspace of $X$ embeds isomorphically into $T$.
\end{description}
\end{example}

Before verifying these assertions we first require some observations.

The norm of $x\in X$, if not equal to $\|x\|_\infty$, is computed by a
tree of sets, the {\em first level\/} being $(E_i)_1^\ell$ where for
some $j$, $(E_i)_1^\ell$ is $j$-admissible and
$$\|x\| = \frac{c}{2^j} \sum_{i=1}^\ell \|E_ix\|\ .$$
For each $i$, if $\|E_ix\|$ does not equal $\|E_ix\|_\infty$, then we
split $\|E_ix\|$ into a second level of sets $m_i$-admissible for some
$m_i$, and so on.  If every set keeps splitting then after $k$ steps we
obtain an expression of the form
\begin{equation}\label{star.two}
  c^k \sum_{s=1}^r 2^{-n(s)} \|F_s x\|\ .
\end{equation}
Of course some sets may stop splitting, in which case if we carry on for
$k$-steps,  we only obtain a lower estimate for $\|x\|$.  Consider the case
where $(x_i)\prec X$ and $x\in \langle x_i\rangle$.  We set $\|x\|_{\T_k,
  (x_i)}$ to be the largest of the expressions of the form
(\ref{star.two}) obtained by splitting $k$-times (a $k$ level tree of
sets, where $(F_s)_1^r$ is the $k^{th}$-level), subject to the additional
constraint that for all $i$ and $s$, $F_s$ does not split $x_i$.  Thus
$F_sx_i$ is either $x_i$ or $0$.
\begin{lemma}\label{E.j}
  Let $(x_i)\prec X$, $\ep>0$ and $k\in\nat$.  Then there exists $x\in
  \langle x_i\rangle$ with  $\|x\|=1$ such that $\|x\|_{\T_k,(x_i)} >1-\ep$.
\end{lemma}
\begin{proof}
  Assume without loss of generality that $\|x_i\| = 1$ for $i\in\nat$.  We
  call $x\in [ x_i]$ an {\em $(n,\ep)$-normalized\/} average (of $(x_i)$
  w.r.t. $(e_i)$) if $x=\sum_{i\in F} a_i x_i/ \|\sum_{i\in F} a_ix_i\|$,
  where $\sum_{i\in F} a_ix_i$ is an $(n,n-1,c\ep/2^n)$-average of
  $(x_i)$ w.r.t.  $(e_i)$.  Thus $(x_i)_{i\in F}$ is $n$-admissible
  w.r.t. $(e_i)$ and if $G\subseteq F$ satisfies $(x_i)_{i\in G}$ is
  $(n-1)$-admissible then $\sum_G a_i < c\ep/2^n$.  Also $\sum_{i\in F}
  a_i=1$ and $a_i>0$ for $i\in F$.  (We can always find such vectors by
  Proposition~\ref{C.l}.)  Note that if $(x_i)_{i \in G}$ is
  $(n-1)$-admissible and if we write $x$ in the form $x=\sum_{i \in F}
  b_ix_i$ (for some $b_i >0$), then $\sum_G b_i < ({c\ep }/{2^n})
  ({2^n}/c) = \ep$ (since $\|\sum_{i\in F} a_i x_i\| \ge c/2^n)$.
 
  We first indicate how to find $x$ satisfying $\|x\|=1$ and
  $\|x\|_{\T_1,(x_i)} >1-\ep$.  Let $\ep_i = 2^{-(i+1)}\ep$ so that
  $\sum_1^\infty \ep_i = \ep/2$.  Let
  $$\|\cdot\|_n =
  \sup \Bigl\{ c 2^{-n} \sum_{j=1}^\ell \|E_j x\| : (E_j x)_{j=1}^\ell
  \mbox{\rm \ is } n \mbox{\rm -admissible} \Bigr\}\ . $$
  and observe that for all $x$,
  $\lim_n\|x\|_n = 0$. Let $n_1 = 1$ and choose $(y_i^1) \prec (x_i)$ and
  $n_j \uparrow \infty$ by induction so that each $y_j^1$ is an
  $(n_j,\ep_j)$-normalized average of $(x_i)$ and for all $j$,
  $\|\sum_{i=1}^j y_i^1\|_{m} <\ep_{j+1}$ if $m \ge n_{j+1}$. 
  Then we choose $y^2$ to
  be an $(n,\ep/2)$-normalized average of $(y_i^1)$ where $n\in\nat$ is
  not important but we may assume that $y^2 = \sum_{F} b_iy_i^1$ where
  $n<n_{\min F}$.
 
  We have $1= \|y^2\|$ and so by the definition of the norm in $X$, there
  exists $j$ such that $1 = \|y^2\|_j = c/2^j \sum_{s=1}^\ell \|E_s(y^2)\|$
  where $(E_sy^2)_1^\ell$ is $j$-admissible.  We claim that by somewhat
  altering the $E_s$'s we can  ensure, by losing no more than $\ep$, that
  the sets $E_s$ do not split any of the $x_i$'s. Indeed if $1\le
  j<n$, then $G= \{i\in F: E_s$ splits $y_i^1$ for some $s\} \in \S_j$.
  Since $j<n$, $\sum_{s\in G} b_s<\ep/2$ and thus by shrinking the
  offending sets $E_s$ to avoid  splitting  $y_i$'s we obtain the desired
  sets.  If $n\le j<n_{\min F}$ then if we fix $i\in F$ and consider
  $G_i=\{r :E_s$ splits one or more of the $x_r$'s in the support of
  $y_i\}$ we get that, by similarly shrinking the offending $E_s$'s so as
  to not split such an $x_r$, and letting $\widetilde E_s$ be the new
  sets,  that
  $$\frac{c}{2^j} \sum \|\widetilde E_s y^2\| > 1-\sum_{i\in F} b_i
  \ep_i > 1-\ep\ .$$
  Finally if $F= (k_1,\ldots, k_r)$ and $n_{k_p} \le j< n_{k_{p+1}}$ then
  $$\Big\| \sum_{\scriptstyle i\in F\atop \scriptstyle i<k_p}
  b_iy_i\Big\|_j < \ep_{k_p} \ \mbox{ and }\ b_{k_p} <\ep/2$$
  so we first discard the $E_s$'s which intersect $\supp (\sum_{i\le
    k_p} b_i y_i)$.  Then arguing as above we shrink the remaining
  $E_s$'s so as to not split any $x_i$.  We obtain
  $$\frac{c}{2^j} \sum \|\widetilde E_s y^2\| \ge 1-\ep_{k_p} - \ep/2 -
  \sum_{\scriptstyle i\in F\atop \scriptstyle i>k_p} b_i
  \ep_i > 1-\ep\ .$$
 
   This proves the lemma in the case $k=1$.  For the general case we
  continue as above letting $(y_i^2)$ be $(n_i^2,\ep_i)$-normalized
  averages of $(y_i^1)$, etc.  If $x=y_1^{k+1}$ then $x$ satisfies the
  lemma for $k$.  We omit the tedious calculations.
\end{proof}

\medskip \pf{Proof of the assertions in Example~\ref{E.h}}
  By Proposition~\ref{D.n}, since $\delta_n(X) \ge c2^{-n}$, we have 
  $\ddot\delta_1 (X) \ge 2^{-1}$.  If there exists
  $\gamma\in \ddot\Delta (X)$ with $\gamma_1\ge 2^{-1}$ then by
  Theorem~\ref{E.gb} there exists $(x_i)\prec (e_i)$ and $d > 0$ so
  that for all $(y_j)_1^\ell \prec (x_i)$ if $m_j = \min (\supp
  (y_j))$ w.r.t.  $x_i$, then
\begin{equation}\label{star.three}
  \Big\|\sum_1^\ell y_j\Big\| \ge d\Big\| \sum_1^\ell \|y_j\| e_{m_j}
  \Big\|_T\ .
\end{equation}
 
  Fix an arbitrary $k$.  By Lemma \ref{E.j} there exists $x\in \langle
  x_i\rangle$ with $\|x\|=1$ and $\|x\|_{\T_k,(x_i)} >1/2$.  Thus there
  exists a $k$-level tree of sets whose final level is $(E_1,\ldots,E_r)$
  so that $c^k \sum_{s=1}^r 2^{-n(s)} \|E_s x\| >1/2$.  Following the same
  partition scheme in $T$ and using (\ref{star.three}) for $y_s = E_s x$
  we get (with   $m_s=\min (\supp (E_sx))$),
  $$d^{-1} = d^{-1}\|x\| \ge \Big\|\sum_{s=1}^r \|E_sx\|e_{m_s} \Big\|_T
  \ge \sum_{s=1}^r 2^{-n(s)} \|E_sx\| > \frac12 (c^{-k})\ .$$ Since
  $c<1$, this is impossible for large enough $k$.  This proves b) for
  $n=1$ and that $\ddot\delta_1(X) =2^{-1}$.  Then Proposition~\ref{D.k}
  yields $\ddot\delta_n(X)=2^{-n}$ for all $n$. 
 
  The remainder of b) easily follows from the proof of
  Proposition~\ref{D.n}. Indeed assume that some $\gamma \in \ddot\Delta(X)$
  satisfies $\gamma_n = 2^{-n}$,   for some $n >1$. By
  Proposition~\ref{D.na}   there is  $(y_i)\prec X $ 
  and an equivalent bimonotone
  norm $|\cdot|$ on $[y_i]$  such that  $(y_i)$ \Dstabs\/
  $\bar\gamma\in \ddot\Delta(X, |\cdot|)$ and $\bar\gamma_n = 2^{-n}$.
  By passing to a
  subsequence we may assume that for some sequence $\ep_n \downarrow
  0$, for all $m$,
  $$\Bigl|\sum_1^k x_i\Bigr| \ge 2^{-n}(1-\ep_m)\sum_1^k |x_i|$$
  if $(x_i)_1^k \prec (y_i)_m^\infty$ and $(x_i)$ is $n$-admissible
  w.r.t. $(y_i)_1^\infty$. Let $\norm{\cdot}$ be the norm constructed
  in the proof of Proposition~\ref{D.n} for $\alpha =1$ and $\theta =
  1/2$.  If $y_r\le x_1 < \cdots < x_r $ then
  $$\Bigl|\sum_1^r x_s\Bigr| \ge (1/2)(1-\ep_r)\sum_1^r |x_s|_{n-1}\
  .$$
  The remaining estimates remain true and, as in the proof of
  Proposition~\ref{D.n},  we obtain
  $$\norm{\sum_1^r x_s} \ge (1/2)(1-\ep_r)\sum_1^r \norm{x_s}\ .$$
  Thus $\gamma_1 = 1/2$  which is impossible.

  If c) were not true, then, by Theorem~\ref{E.a} b), a subspace $Y$
  of $X$ isomorphic to a subspace of $T$ would admit a renorming for
  which $\gamma_1 (Y) = 1/2$, in contradiction to b).
\endpf

\begin{remark}\label{E.k}\rm
  The above example $X$ yields the following.  There exists
  $(x_i)\prec (e_i)$ and a sequence of equivalent norms
  $\norm{\cdot}_j$ so that for all $k$ on $[x_i]_k^\infty$,
  $\|x\| \ge \norm{x}_j \ge c^2\|x\|$ if $j\ge
  k$ and furthermore $\delta_1 (\norm{\cdot}_j, (x_i)) >\frac12
  -\ep_j$ for some $\ep_j\to0$.  Yet $\gamma_1 <\frac12$ for all
  $\gamma\in \ddot\Delta (X)$.  To see this one needs only choose
  $(x_i)$ so that on $[x_i]_k^\infty$, $\|x\|
  =\sup_{\ell\ge k} \|x\|_\ell$.  This can be accomplished by taking
  each $x_j$ to be an iterated $j+1$-normalized average of $(e_i)$ (as
  in lemma~\ref{E.j}).  Then set $\norm{x}_j = (1/j) \sum_1^j
  \|x\|_i$.  Since $\|x\| \ge \|x\|_i \ge c\|x\|_j \ge c^2 \|x\|$ on
  $\langle x_s\rangle_j^\infty$, $\|x\| \ge \norm{x}_j \ge c^2 \|x\|$.
\end{remark}

We mention one other example, taken from \cite{AO}.
First suppose that $X= T(\S_n,\ov{\theta}_n)_{n\in\nat}$ where
$1>\sup_n \ov{\theta}_n$ and $\lim_{n\to\infty} \ov{\theta}_n =0$.  We
shall call $(\theta_n)$ {\em regular\/} if for all $n,m\in \nat$,
$\theta_{n+m} \ge \theta_n \theta_m$.  It is easy to verify that every
such $X$ has a regular representation, i.e., for some regular sequence
$(\theta_n)$ we have $X= T(\S_n,\theta_n)_{\nat}$. Thus $\lim_n
\theta_n^{1/n}$ exists by Lemma~\ref{D.l}.

\begin{example}\label{F.e}
  Let $X = T(\S_n,\theta_n)_{\nat}$ where $1>\sup_n\theta_n$,
  $\theta_n\to 0$ and $(\theta_n)$ is regular.  Let $\theta = \lim_n
  \theta_n^{1/n}$.  Then
  \begin{description}
  \item[a)] For all $Y \prec X$  we have $\ddot{\delta}_1(Y)=\theta$.
  \item[b)] For all $Y \prec X$ and for all $n \in \nat$,
    $\ddot{\delta}_n(Y)=\theta^n$     and $\ddot{\delta}_{\omega}=0$.
  \item[c)] For all  $Y \prec X$, $ I_{\Delta}(Y) =
     \left\{ \begin{array}{ll}
      \omega & \mbox{ if } \theta =1 \\
      1 & \mbox{ if } \theta < 1 
      \end{array} \right. $
  \item[d)] For all $Y \prec X$ and $j \in \nat$ we have $\delta_j(Y)
    \le \theta^j \sup_{n \ge j} \theta_n \theta^{-n} \vee
    \theta_j/\theta_1$.
    In particular,
    if ${\theta_n}{\theta^{-n}} \to  0$ then $X$ is
    arbitrarily distortable.
\end{description}
\end{example}

\section{Renormings of  $T$,  and spaces of bounded distortion}

\begin{defin}\label{F.b}
  The distortion constant of a space $X$ is defined by
  $$D(X) = \sup_{|\cdot|\sim \|\cdot\|} d(X,|\cdot|)\ . $$
\end{defin}

So $X$ is distortable iff $D(X)>1$. Similarly, $X$ is arbitrarily
distortable iff $D(X) = \infty$. Finally, $X$ is of bounded
distortion iff there is $D< \infty$ such that $D(Y) \le D$ for every
subspace $Y \subseteq X$.

As we saw in Proposition~\ref{B.d}, Tsirelson's space $T$ satisfies $D(T)
\ge 2$.  Similarly one can show that $D(T_\theta)\ge \theta^{-1}$.  However,
not much more is  known about distorting  $T$. 
It is unknown if $T$ is arbitrarily distortable, or at least whether
it contains an arbitrarily distortable subspace; and, if  not,
what  is  $D(T)$ or at least a reasonable upper
estimate for it.  The interest in these questions lies in the fact
that, as already mentioned, no examples are yet known of distortable
spaces which are of  bounded distortion.

From techniques developed earlier in this paper we easily get some
information on asymptotic constants of equivalent norms on Tsirelson
space. This should be compared with Theorem~\ref{E.a} where the constants
for the original norm were established.

Surprisingly, it is not known if there exists $(x_i)\prec T$ and
an equivalent norm $|\cdot|$ on $[x_i]$ with
$\delta_1((x_i),|\cdot|)<1/2$.  Our next result shows that the class of
equivalent norms for which $\delta_1 =1/2$ cannot arbitrarily distort
$T$.

\begin{thm}\label{E.b}
  There exists an absolute constant $D$ with the following property.
  Let $X\prec T$ and let $|\cdot|$ be an equivalent norm on $X$ such
  that for some $\gamma \in\Delta (X,|\cdot|)$, $\gamma_1=1/2$.  Then
  $d(X,|\cdot|)\le D$.
\end{thm}
\begin{proof}
  Let $(z_i)$ be a basic sequence in $X$ $\Delta$-stabilizing $\gamma$
  under $|\cdot|$ where $\gamma_1=1/2$.  Let $\ep >0$.  By passing to
  a block basis of $(z_i)$ and multiplying $|\cdot|$ by a constant if
  necessary we may assume that $\|\cdot\|_T \ge |\cdot|$ on $[z_i]$
  and for all $(w_i)\prec (z_i)$ there exists $w\in \langle
  w_i\rangle$ with $1+\ep > \|w\|_T \ge |w|=1$.  Choose a normalized
  block basis $(w_i)$ of $(z_i)$ satisfying $1+\ep \ge \|w_i\|_T \ge
  |w_i| =1$ for all $i$.  Theorem~\ref{E.gb} allows us to also assume
  that
  $$| \sum a_iw_i| \ge  (1/2-\ep) \|\sum a_i e_i\|_T\ .$$
  There exists an
  absolute constant $D_1$ so that $(w_i/\|w_i\|_T)$ is $D_1$-equivalent
  to $(e_{m_i})$ in $\|\cdot\|_T$, where $m_i = \min\supp (w_i)$ w.r.t.
  $(e_i)$, for each $i$ \cite{CJT}.  Thus we have, for all $(a_i)\subseteq
  \real$,
\begin{equation}\label{star}
  (1+\ep) D_1 \|\sum a_i e_{m_i}\|_T \ge \|\sum a_iw_i\|_T \ge |\sum
  a_iw_i| \ge (1/2-\ep) \|\sum a_i e_i\|_T\ .
\end{equation}
 
  Consider the subsequence $(p_i)$ of $\nat$ defined by induction by $p_1
  = 1$ and $p_{i+1} = m_{p_i}$, for $i \ge 1$.  There is a universal
  constant $D_2$ so that $(e_{p_i})$ is $D_2$-equivalent to
  $(e_{p_{i+1}})$ in $\|\cdot\|_T$ \cite{CJT}.  Also, on the  subspace
  $[w_{p_i}]$ we have, by (\ref{star}),
  $$(1+\ep)D_1\|\sum a_i e_{p_{i+1}}\|_T \ge |\sum a_i w_{p_i}|
  \ge  (1/2-\ep)   \|\sum a_i e_{p_i}\|_T\ .$$
  Thus the conclusion follows with $D= 2D_1D_2$.
\end{proof}

A natural question in light of the above results is whether one can
quantify the distortion $d(X,|\cdot|)$ of an equivalent norm $|\cdot|$
on $X\prec T$ in terms of $\Delta (X,|\cdot|)$.

\begin{problem}\label{E.d}
  Let $|\cdot|$ be an equivalent norm on $T$ and let $(x_i)\prec T$
  $(\Delta,|\cdot|)$-stabilize $\gamma$.  Thus for some $c >0$, $c2^{-n}
  \le \gamma_n \le 2^{-n}$ for all $n$.  Does there exist a function
  $f(c)$ so that $d(X,|\cdot|) \le f(c)$?
\end{problem}

We shall give a suggestive partial answer to a weaker problem.  First
we note the following proposition.

\begin{prop}\label{E.da}
  For $n \in \nat$ define the equivalent norm $\|\cdot\|_n$ on $T$ by
  $\|x\|_n = \sup \{ 2^{-n} \sum_1^\ell \|E_i x\| :(E_i x)_1^\ell$ is
  $n$-admissible$\}$. Given $X \prec T$ and $\ep_n \downarrow 0$ there
  exists $(x_i) \prec X$ so that for all $n$ if $x \in \langle
  x_i\rangle_n ^\infty$  then $\Bigl|\|x\| - \|x\|_n \Bigl| < \ep_n \|x\|$.
\end{prop}

\begin{proof}
  First note that if
  $$\|\cdot\|_{\S_n} = \sup \Bigl\{\sum_{i\in E} |x(i)| :
   E\in \S_n\Bigr\}$$
  then for all $x \in T$ we have $\|x\|_n \le \|x\|
  \le \|x\|_n + \|x\|_{\S_n}$. Indeed, if $\|x\| \ne \|x\|_\infty$
  then $\|x\| = x^*(x)$ for some functional $x^*$ (with $\|x^*\| = 1$)
  determined by the successive iterations of the implicit equation of
  the norm in $T$; in particular, $x^* (e_i) = \pm 2^{- n(i)}$ for all
  $i$. We may write $x^* = y^* + z^*$ where $z^*(e_i) = \pm 2^{-
    n(i)}$ if $n(i) \le n$ and 0 otherwise. Thus, since the support of
  $z^*$ is $n$-admissible, $|z^* (x)| \le (1/2)\|x\|_{\S_n}$ and
  $|y^*(x)| \le \|x\|_n$. Furthermore, $\|x\|_{\S_n} \le 2^n \|x\|$.
  Since the Schreier space $S_n$ is isomorphic to a subspace of
  $C(\omega^{\omega^n})$ (Remark~\ref{C.e}),
, it is $c_0$-saturated, i.e.,  every
  infinite-dimensional subspace contains a copy of $c_0$, and thus
  $\|\cdot\|_{\S_n}$  cannot be equivalent to $\|\cdot\|$  on any
  infinite-dimensional subspace of $T$. In particular we can chose
  $(x_i)\prec X$ so that for all $x \in \langle x_i\rangle_n ^\infty$,
  $\|x\|_{\S_n} \le \ep_n \|x\|$. The conclusion follows.
\end{proof}

\begin{problem}\label{E.e}
  Let $|\cdot|$ be an equivalent norm on $X= [x_i]\prec T$.  Let
  $(y_i)\prec (x_i)$, $C<\infty$ and suppose that for all $n$, if $y\in
  [y_i]_n^\infty$ then $C^{-1}|y|_n\le |y|\le C|y|_n$,
  where $|y|_n = \sup \{2^{-n} \sum_1^\ell |E_iy| :(E_i y)_1^\ell$ is
  $n$-admissible w.r.t. $(x_i)\}$.  Does there exist a function $F(C)$
  so that $d(Y,|\cdot|)\le F(C)$?
\end{problem}
\begin{prop}\label{E.f}
  Let $(y_i)\prec (x_i)\prec T$ and let $|\cdot|$ be an equivalent norm
  on $[x_i]$.  Suppose that for all $n$ and $y\in [y_i]_n^\infty$,
  $C^{-1}|y|_n \le |y| \le C|y|_n$ (where $|\cdot|_n$ is  defined as
  above).  Then for all $\ep > 0$ there exists $n_0$ and an equivalent
  norm $\norm{\cdot}$ on $[ y_i]_{n_0}^\infty$ such that $C^{-1}\norm{y}
  \le |y| \le C\norm{y}$ for $y \in [y_i]_{n_0}^\infty$ and
  $\delta_1 ((y_i)_{n_0}^\infty,\norm{\cdot}) >\frac12 -\ep$.
\end{prop}

\begin{proof}
  Choose $n_0$ so that ${C^2}/n_0 <\ep$.  On $[y_i]_{n_0}^\infty$
  define $\norm{y} = \frac1{n_0} \sum_1^{n_0} |y|_j$.  Clearly the
  inequality between the norms hold.  Let $p\in\nat$ and let
  $(z_i)_1^p \prec [y_i]_{n_0}^\infty$ satisfy $y_{n_0+p} \le z_1
  <\cdot < z_p$.  Let $z= \sum_1^p z_i$.  Then (see the proof of
  Proposition~\ref{D.n}) $|z|_{j+1} \ge \frac12 \sum_{i=1}^p |z_i|_j$
  for $j=1,\ldots,n_0-1$.  Hence
  $$   \norm{z}  \ge  \frac1{n_0} \sum_{j=1}^{n_0-1} \frac12
  \sum_{i=1}^p |z_i|_j
  = \frac12 \sum_{i=1}^p \norm{z_i} - \frac1{2{n_0}} \sum_{i=1}^p
  |z_i|_{n_0}\ .$$
  Now $|z_i|_{n_0} \le C|z_i| \le C^2 \norm{z_i}$ and so
  $$  \norm{z}  \ge  \frac12 \sum_{i=1}^p \norm{z_i}
  \left(1-\frac{C^2}{2{n_0}}
  \right) >  (1-\ep) \frac12 \sum_{i=1}^p \norm{z_i} \ ,$$
  completing  the proof.
\end{proof}

Finally, let us recall the following known \cite{CJT} property of $T$.
There exists an absolute constant $D_1$ so that if $x_1<y_1<x_2<y_2<\cdots$
are normalized in $T$ then $(x_i)$ is $D_1$-equivalent to $(y_i)$.
It turns out that equivalent norms on $T$ that satisfy this property
(with a fixed  constant) cannot arbitrarily distort $T$.
The result, in fact, holds in any space having this subsequence
property.

\begin{prop}\label{E.l}
  There exists a function $f(D)$ satisfying the following.  If $|\cdot|$
  is an equivalent norm on $[x_i]_{\nat}\prec T$ so that $(y_i)$ is
  $D$-equivalent to $(z_i)$ whenever $y_1<z_1<y_2<\cdots$ is a normalized
  block basis of $(x_i)$, then $d(X,|\cdot|) \le f(D)$.
\end{prop}
\begin{proof}
  By passing to a block basis of $(x_i)$ and scaling the norm $|\cdot|$
  we may assume that there exists $d > 1$ so that for all $x\in [x_i]$,
  $d^{-1} \|x\| \le |x| \le \|x\|$; furthermore, in any block subspace $Y$
  of $(x_i)$ there exist $y,z \in Y$ with $|y|=|z|=1$ and $\|y\| \le 2$
  and $\|z\|>d /2$.  Choose a $|\cdot|$-normalized block basis of
  $(x_i)$, $y_1<z_1<y_2<\cdots$ with $\|z_i\| > d/2$ and $\|y_i\|\le 2$
  for all $i$.  There exists $w = \sum a_i z_i$ satisfying $|w| =1$ and
  $\|w\|<2$.  Since $(z_i)$ and $(y_i)$ are $D$-equivalent for $|\cdot|$,
  $|\sum a_i y_i| > D^{-1}$.  Also $(z_i/\|z_i\|_T)$ and
  $(y_i/\|y_i\|_T)$ are $D_1$-equivalent in $T$.  Thus
  $$  \|\sum a_i y_i\|_T  \le  2D_1 \|\sum a_i z_i /\|z_i\|_T\,\|_T
  \le  4D_1 /d\ \| \sum a_i z_i\|_T \le 8D_1/ d\ .$$
  Thus $D^{-1} \le 8D_1/d$ and so $d\le 8D_1 D \equiv f(D)$.
\end{proof}

We now turn to some results about spaces of bounded distortion.

\begin{thm}\label{F.f}
  Let $X$ be an asymptotic $\ell_1$ space.  Let $\gamma\in \Delta (X)$
  and let $(y_i)\prec X$ \Dstab\/ $\gamma$.  If
  $Y=[y_i]$ is of $D$-bounded distortion then for any $\alpha
  <\omega_1$ and $n,m\in\nat$,
\begin{description}
\item[a)] $D^{-1} (\ddot\delta_\alpha (Y))^n \le \gamma_{\alpha\cdot
    n} \le (\ddot\delta_\alpha (Y))^n$
\item[b)] $\gamma_{\alpha\cdot n} \gamma_{\alpha\cdot m} \le
  \gamma_{\alpha\cdot(n+m)} \le D^2 \gamma_{\alpha\cdot n}
  \gamma_{\alpha\cdot m}$.
\end{description}
\end{thm}
\begin{proof}
  a) Let $\ov{\gamma} = (\,\ov{\gamma}_\alpha) \in \ddot\Delta (Y)$.
  Choose an equivalent norm $|\cdot|$ on $Y$ and $(w_i) \prec (y_i)$
  which $(\Delta, |\cdot|)$-stabilizes $\ov{\gamma}$.  Let $\ep>0$.
  By passing to a block basis of $(w_i)$ and scaling $|\cdot|$ we may
  suppose that
  $$|w| \le \|w\| \le (D+\ep) |w|\ \mbox{\rm  for } w\in [w_i]\ .$$

  Let $\alpha < \omega_1$ and $n \in \nat$.
  We  may  assume that $\delta_{\alpha\cdot n}
  ((w_i),|\cdot|) > \ov{\gamma}_{\alpha\cdot n} - \ep$.  Thus if
  $(x_s)_1^r$ is $\alpha\cdot n$-admissible w.r.t. $(w_i)$,
  $$ \Big\| \sum_1^r x_s\Big\|  \ge  \Big| \sum_1^r x_s\Big| \ge
  (\,\ov{\gamma}_{\alpha\cdot n}-\ep)
  \sum_1^r |x_s|
  \ge  {\ov{\gamma}_{\alpha\cdot n}-\ep \over D+\ep} \sum_1^r
  \|x_s\|\ .$$
 
  It follows that $\gamma_{\alpha\cdot n}\ge
  \ov{\gamma}_{\alpha\cdot n} /D$ and so $ \ov{\gamma}_{\alpha\cdot n}\le
  D\gamma_{\alpha\cdot n}$.  Passing to the supremum over all $
  \ov{\gamma}_{\alpha\cdot n}$ and using Proposition~\ref{D.j} g), we get
  $\ddot\delta_{\alpha\cdot n} (Y) \le D\gamma_{\alpha\cdot n}$.  Hence
  by Proposition~\ref{D.k},
  $$   D^{-1} (\ddot\delta_\alpha (Y))^n = D^{-1} \ddot\delta_{\alpha
    \cdot  n} (Y)
  \le \gamma_{\alpha\cdot n}  \le  \ddot\delta_{\alpha\cdot n} (Y)
  = (\ddot\delta_\alpha (Y))^n\ . $$

b) Using part a) and Proposition~\ref{D.j} d),
\begin{eqnarray*}
  \gamma_{\alpha\cdot n} \gamma_{\alpha\cdot m} & \le &
  \gamma_{\alpha\cdot (n+m)} \le (\ddot\delta_\alpha (Y))^{n+m}\\
  & = & (\ddot\delta_\alpha (Y))^n (\ddot\delta_\alpha (Y))^m \le D^2
  \gamma_{\alpha\cdot n} \gamma_{\alpha\cdot m}\ ,
\end{eqnarray*}
completing the proof.
\end{proof}

Combining the proposition with Theorem~\ref{F.c} we get a complete
description, up to equivalence, of sequences $\gamma$ from $ \Delta(X)$,
in spaces of $D$-bounded distortion. We leave the details to the reader.

Recall the notation $\widehat\gamma_\alpha = \lim_k
(\gamma_{\alpha\cdot k})^{1/k}$, for $\alpha < \omega_1$
(Corollary~\ref{D.la}).  If $Y \prec X$ \Dstabs\/ $\gamma$,
we may write $\widehat\gamma_\alpha (Y)$ to emphasise the subspace
$Y$.  By Proposition~\ref{D.m}, $\widehat\gamma_\alpha (Y)=
\ddot\delta_\alpha (Y)$. Therefore, by Proposition~\ref{A.d}, we have
an important sufficient condition for an asymptotic $\ell_1$ space to
contain an arbitrary distortable subspace.
\begin{cor}\label{F.fa}
  Let $X$ be an asymptotic $\ell_1$ space.  Let $\gamma\in \Delta (X)$
  and let $(y_i)\prec X$ $\Delta$-stabilize $\gamma$.  If there exists
  $\alpha < \omega_1$ such that $\gamma_\alpha >0$ and $\lim_n
  \gamma_{\alpha \cdot n} \widehat\gamma_\alpha (Y)^{-n} =0$, then $Y$
  contains an arbitrarily distortable subspace.
\end{cor}

Let us present an alternative approach to Corollary~\ref{F.fa}, taken
from \cite{Tom}, which is of independent interest.  It is based on a
construction of certain asymptotic sets in a general asymptotic
$\ell_1$ space.

\medskip \pf{An alternative proof of Corollary~\ref{F.fa}} {\it
  (Sketch)\ \ }
  Let $\gamma\in \Delta (X)$, let $Y= [y_i]\prec X$
  $\Delta$-stabilize $\gamma$ and  let $(y_i^*)$ be the biorthogonal
  functionals in $Y^*$.  Suppose that $Y$ is of $D$-bounded distortion.
  Fix an  arbitrary $\alpha < \omega_1$.  We shall
  show that $(1/3D) (\widehat\gamma_\alpha(Y))^n \le \gamma_{\alpha\cdot
  (n-1)}$.  By Proposition~\ref{D.m}, this is slightly weaker than
  Theorem~\ref{F.f}, but sufficient to imply Corollary~\ref{F.fa}.
 
  Fix $n \in \nat$.  First we shall show that for all $\ep >0$, all
  normalized blocks $(x_i)\prec (y_i)$, and all $ 0<\lambda < 1$,
  there is an
  $(\alpha\cdot n, \alpha \cdot (n-1), \ep)$ average $x$ of $(x_i)$
  w.r.t. $(y_i)$ such that $\|x\| \ge \lambda (\widehat\gamma_\alpha(Y))^n
  \equiv \lambda'$.
 
  This is done by blocking, in the spirit of James \cite{Jam}.
  Fix $m$ sufficiently large and pick
  $N \subseteq \nat$ such that $[\S_{\alpha\cdot n}]^m (N)\subseteq
  \S_{\alpha\cdot (n\,m)}$ (Corollary~\ref{C.d})  and that
  $\lambda \gamma_{\alpha \cdot(n\,m)} \le \delta_{\alpha
      \cdot(n\,m)}((x_i)_N)$ (this is possible by the
  Definition~\ref{D.g} of the $\Delta$-spectrum).
  Pick
  $(z_i^{(1)})\prec (x_i)_N$ such that for all $i$, $z_i^{(1)}$ is an
  $(\alpha\cdot n, \alpha \cdot (n-1), \ep)$ average of $(x_i)_N$
  w.r.t.  $(y_i)$. If for all $i$, $\|z_i^{(1)}\| < \lambda'$, then pick
  $(z_i^{(2)})\prec (z_i^{(1)})$ such that for all $i$, $z_i^{(2)}$ is
  an $(\alpha\cdot n, \alpha \cdot (n-1), \ep)$ average of
  $(z_i^{(1)}/\|z_i^{(1)}\|)$ w.r.t. $(y_i)$.  And keep going. Assume
  that after $m$ steps we still had that $\|z_i^{(k)}\| < \lambda'$ for
  all $i$ and all $k \le m$. Write $z_1^{(m)}= \sum_{j \in N} b_j
  x_j$; then $b_j \ge 0$ and let $J$ be the set of all $j\in N$ such
  that $b_j > 0$. It is easily seen that $(x_j)_{j \in J}$ is
  $[\S_{\alpha\cdot n}]^m (N)$-admissible w.r.t. $(y_i)$, hence also
  $(\alpha\cdot(n\,m))$--admissible w.r.t. $(y_i)$. Moreover, our
  assumption on the norms of the $z_i^{(k)}$'s easily yields that
  $\sum b_j > (1/\lambda')^{m-1}$. Thus
  \begin{eqnarray*}
    \lefteqn{
      (1/\lambda')^{m-1} \lambda \gamma_{\alpha\cdot(n\,m)} \le
  \lambda \gamma_{\alpha\cdot(n\,m)} \sum_{j \in J} \|b_j x_j\|}\\
  && \qquad \le  \delta_{\alpha  \cdot(n\,m)}((x_i)_N)   
  \sum_{j \in J} \|b_j x_j\|
  \le \|\sum_{j \in J} b_j   x_j\| = \|z_1^{(m)}\| < \lambda'\ .  
  \end{eqnarray*}
  It follows that $\lambda\gamma_{\alpha\cdot(n\,m)} < \lambda'^m$, hence
  $(\gamma_{\alpha\cdot(n\,m)} )^{1/n\,m} < \lambda^{1/n - 1/m n}
  \widehat\gamma_\alpha(Y)$, a contradiction, if $m$ is large enough.
 
  Now we shall define asymptotic sets $A, B \subseteq S(Y)$ and a set
  $A^*$ in the unit ball of $Y^*$ such that $A^*$ 2-norms $A$ and the
  action of $A^*$ on $B$ is small.  By passing to a tail subspace of
  $Y$ if necessary, we may assume without loss of generality that
  $\frac34 \gamma_{\alpha \cdot (n-1)} \le \frac78 \delta_{\alpha
    \cdot (n-1)}(Y)$.  Fix $\ep >0$, quite small as determined at the
  end of this  proof.  Let $A^*$ consist of all functionals in $Y^*$ of
  the form $y^* = \frac34 \gamma_{\alpha \cdot (n-1)}\sum_{k \in K}
  w_k^*$, where $(w_k^*)_K \prec (y_i^*)$ is $(\alpha \cdot
  (n-1))$-admissible (w.r.t. $(y_i^*)$); and let $A$ consist of all $y
  \in S(Y)$ that are $2$-normed by $A^*$. The set $A$ is asymptotic by
  the definition of the $\Delta$-stabilization.  Since $Y$
  \Dstabs\/ $\gamma$, it is not difficult to see that $A$ is
  asymptotic in $Y$ and that functionals from $A^*$ have the norm not
  exceeding 1. Then $B$ consists of all vectors of the form $x/\|x\|$,
  where $x$ is an $(\alpha\cdot n, \alpha \cdot (n-1), \ep)$ average
  w.r.t. $(y_i)$ of some normalized $(x_i)\prec (y_i)$, such that
  $\|x\| \ge (1 - \ep) (\widehat\gamma_\alpha(Y))^n $.  By the first
  part of this proof, $B$ is asymptotic in $Y$. We will show that if
  $y^* \in A^*$ and $z \in B$, then $|y^* (z)|\le \frac34
  \widehat\gamma_\alpha(Y)^{-n}
  (\gamma_{\alpha \cdot (n-1)}+ \frac76 \ep)/ (1 -   \ep) \equiv \eta$.
 
  This is a direct consequence of the following estimate.
  If $x$ is an $(\alpha\cdot n, \alpha \cdot (n-1),
  \ep)$ average as above, and if $(E_k) \in \S_{\alpha \cdot (n-1)}$,
  and $E_k x$ denotes the restriction of $x$ whose support w.r.t.
  $(y_i)$ is $E_k$; then $\sum \|E_k x\| \le 1+7\ep /6\gamma_{\alpha
    \cdot (n-1)}$.  To see this, write $x$ in the form
  $x = \sum_{i\in F} a_i x_i$
  where $(x_i)_{i\in F}$ is $\alpha \cdot n$-admissible w.r.t.
  $(y_i)$ and if $J\subseteq F$ satisfies $(x_i)_{J}$ is $\alpha
  \cdot(n-1)$-admissible then $\sum_G a_i < \ep$.  Also $\sum_{i\in F}
  a_i=1$ and $a_i>0$ for $i\in F$.  Set $I= \{i : E_k\cap \supp (x_i)
  \ne \emptyset$ for at most one $k\}$ and $J= F\setminus I$; and for
  $i \in J$ let $K_i= \{k : E_k\cap \supp (x_i) \ne \emptyset\}$. Then
  it can be checked that $(x_i)_J$ is $\alpha \cdot (n-1)$-admissible,
  hence
  $$\sum_k \|E_k x\| \le \sum_{i \in I} a_i \|x_i\| + \sum_{i \in J} a_i
  \sum_{k \in K_i}\|E_k x_i\| \le 1 + \ep/ \delta_{\alpha \cdot
    (n-1)}(Y) \le 1 + 7\ep/ 6\gamma_{\alpha \cdot
    (n-1)}\ .$$
  Now, if $y^* = \frac34 \gamma_{\alpha \cdot (n-1)}\sum_{k \in K}
  w_k^* \in A^*$ then letting $E_k = \supp(w_k)$ for all $k$ we get
  $|y^*(z)| \le \eta $, as required.
 
  As mentioned in Section 2, $Y$ is $(1/2+ 1/4\eta)$-distortable.  Hence
  the assumption of $D$-bounded distortion implies $1/2 + 1/ 4\eta\le D$.
  Substituting the definition of $\eta$ and taking $\ep >0$ sufficiently
  small we get the inequality $(1/3D) (\widehat\gamma_\alpha(Y))^n \le
  \gamma_{\alpha\cdot (n-1)}$, as promised.
\endpf

As we remarked earlier, the assumption of bounded distortion implies
the existence of certain subspaces with a nice structure (\cite{MiT},
\cite{Ma}, \cite{To}).  We would like to identify more such regular
subspaces in the class of asymptotic $\ell_1$ spaces of bounded
distortion.

Recall (Proposition~\ref{E.da}) that in Tsirelson's space $T =
T_\theta$, for all $\ep_n \downarrow 0$ there exists $(x_i) \prec T$
so that for all $n$  and all $x \in \langle x_i\rangle_n^\infty$ we
have
$$ (1+\ep_n)^{-1}\|x\|_T \le
\sup \{ \theta^{n} \sum_1^\ell \|E_i x\|_T :(E_i)_1^\ell
\mbox{\rm  \ is } n \mbox{\rm -admissible}\} \le \|x\|_T\ . $$

In any asymptotic $\ell_1$ space with bounded distortion one can find a
block basis that displays  an isomorphic version of this  phenomenon.

\begin{thm}\label{F.l}
  Let $X$ be an asymptotic $\ell_1$ space of $D$-bounded distortion not
  containing $\ell_1$. There exist $(w_i)\prec X$, $\alpha =
  \omega^{\beta_0}$, $0 < \theta <1$, and   $(z_i) \prec
  (w_i)$ such that for every $k\in \nat$ we have, for $z \in [z_i]_k^\infty$,
\begin{eqnarray*}
\lefteqn{%
  (1/4\,D) \sup_{1\le n \le k} \sup\Bigl\{\theta^n \sum \|E_i z\| :
  (E_i) \mbox{\rm  \ is } \alpha\cdot n \mbox{\rm -admissible} \Bigl\}
  \le \|z\|}\\
  & &\qquad \le 4\, D \inf_{1\le n \le k} \sup\Bigl\{ \theta^n\sum
  \|E_i   z\| : (E_i)  \mbox{\rm  \ is }
  \alpha\cdot n \mbox{\rm -admissible} \Bigl\}\ .
\end{eqnarray*}
(Here, for an interval $E$ of $\nat$ and $z = \sum a_i w_i \in [w_i]$,
$Ez$ denotes the restriction w.r.t. $(w_i)$, i.e., $E z = \sum_{i \in
  E} a_i w_i $.)
\end{thm}
\begin{proof}
  By Proposition~\ref{D.f}, $\ddot\delta_\beta(X) > 0$ for at most
  countably many $\beta$'s; write this set as $(\beta_m)$.  For an
  arbitrary $\beta < \omega_1$, it follows from Lemma~\ref{D.e} that if
  $(y_i) \prec (e_i)$ then $\ddot\delta_\beta ((y_i)_n^\infty) =
  \ddot\delta_\beta (y_i)$ for all $n$; and that $\ddot\delta_\beta
  (z_i)\le \ddot\delta_\beta (y_i)$ whenever $(z_i) \prec (y_i)$.
  Letting, for example, $f(y_i) = \sum 2^{-m} \ddot\delta_{\beta_m}
  (y_i)$, by a standard induction argument, similar to that  in
  Proposition~\ref{D.h}, we can stabilize $f(y_i)$.  That is, we can
  find $(y_i)\prec X$ such that $f(z_i)= f(y_i)$ for all $(z_i) \prec
  (y_i)$.  Since $\ddot\delta_\beta(X) = 0$ implies $\ddot\delta_\beta
  (z_i) =0$ for all $(z_i)\prec X$,  the stabilization of $f$
  implies that we have, for all $(z_i) \prec (y_i)$,
  $$ \ddot\delta_\beta (z_i) = \ddot\delta_\beta (y_i) \qquad\mbox{\rm for
  all  }\beta < \omega_1\ .$$
 
  Let $\alpha = I_\Delta (y_i)$; by Theorem~\ref{F.c}, $\alpha =
  \omega^{\beta_0}$ for some $\beta_0 < \omega_1$. Let
  $\theta = \ddot\delta_\alpha (y_i)$.
  Then $\ddot\delta_{\alpha\cdot n} (y_i) = \theta^n$ for $n \in
  \nat$,   by Proposition~\ref{D.k}.
  By an inductive construction followed by a diagonal argument,
  using Proposition~\ref{D.na}, we can
  find $(w_i) \prec (y_i)$ and equivalent bimonotone norms $|\cdot|_n$ on
  $[w_i]_n^\infty$ such that for all $(z_i) \prec (w_i)_n^\infty$ and $n
  \in \nat$,
  \begin{equation}   \label{st}
    \delta_\alpha ([z_i]_n^\infty, |\cdot|_n)\ge  2^{-1/n}\theta.
  \end{equation}
  Notice that (\ref{st}) is  preserved if the norms
  involved are multiplied by constants. Therefore by scaling and the
  assumption of bounded distortion we may additionally ensure that
  $\|w\|\le |w|_n \le 2\,D\|w\|$ for $w \in [w_i]_n^\infty$ and all
  $n\in \nat$.
   
  Now, given any $\alpha$-admissible family of intervals $(F_i)_1^{k}$
  of $\nat$, let $(G_i)_1^{k}$ be adjacent intervals such that $\min
  F_i = \min G_i$ for $i < k$ and let $G_k = F_k$.  Since the norms
  $|\cdot|_n$ are bimonotone, $|F_iw|_n \le |G_i w|_n$ for $w \in
  [w_i]_n^\infty$ and all $n \in \nat$.  In particular, by
  Remark~\ref{D.ca}  for $n\in\nat$ and $w \in [w_i]_n^\infty$ we get
  $$
     |w|_n \ge \delta_\alpha\Bigl([w_i]_n^\infty, |\cdot|_n\Bigr)
  \sum_{i=1}^k    |G_i w|_n \ge \delta_\alpha\Bigl([w_i]_n^\infty,
  |\cdot|_n\Bigr)   \sum_{i=1}^k |F_i w|_n \ .
  $$

    Using this and the assumption (\ref{st}) on $\delta_\alpha$'s we
  easily get, for $n\in\nat$ and $w \in [w_i]_n^\infty$,
  \begin{eqnarray*}
    2\,D \|w\|\ge |w|_n &\ge & \sup\Bigl\{\delta_\alpha^n \sum_{i=1}^k
       |E_i w|_n : (E_i)  \mbox{\rm  \ is } [\S_\alpha]^ n
       \mbox{\rm -admissible} \Bigl\} \\
    &\ge& (1/2) \sup\Bigl\{\theta^n \sum_{i=1}^k \|E_i w\| :
    (E_i)  \mbox{\rm  \ is } [\S_\alpha]^ n
       \mbox{\rm -admissible}\Bigl\}\ ,
  \end{eqnarray*}
  where we have abbreviated $ \delta_\alpha \Bigl([w_i]_n^\infty,
  |\cdot|_n\Bigr )$ to   $\delta_\alpha $.
  Finally, using  Corollary~\ref{C.d} and a diagonal argument,
  construct  a subsequence   $M = (m_i) $ of $\nat$ such that setting
  $M_n = (m_i)_n^\infty$ we get $\S_{\alpha \cdot n}(M_n)
  \subseteq [\S_\alpha]^n$ for all $n$. Thus for $w \in
  [w_i]_{i\in{M_n}}$ replacing the  supremum in the last formula
  by the supremum over $\S_{\alpha \cdot n} (M_n)$-admissible
  families and relabelling the subsequence by  $(w_i')$  we get,
  for $n \in \nat$  and $w \in [w_i']_n^\infty$,
  $$\|w\| \ge (1/4\,D)  \sup\Bigl\{\theta^n
  \sum \|E_i w\| : (E_i)
  \mbox{\rm  \ is } \alpha\cdot n   \mbox{\rm -admissible}  
  \Bigl\}\ .$$
  It should be noted that in this last estimate, the admissibility condition
  is understood with respect to the above  subsequence $(w_i')$
  of $(w_i)$    which indeed corresponds to the  subsequence $M$ of $\nat$.

  We relabel once more, denoting $(w_i')$ simply by $(w_i)$.
  Set $\norm{w}_n = \sup\Bigl\{\theta^n \sum \|E_i w\| : (E_i)
  \mbox{\rm  \ is } \alpha\cdot n   \mbox{\rm -admissible}\Bigl\}$,
  for $w \in [w_i]_n^\infty$ and
  $n\in\nat$. These are equivalent norms on the subspaces where they are
  defined.  Therefore stabilizing all norms $\norm{\cdot}_n$ on a
  nested sequence of block subspaces, using the assumption of bounded
  distortion, and passing to a diagonal subspace we get $(z_i) \prec
  (w_i)$ and $A_n$ such that $[z_i]_n^\infty \prec [w_i]_n^\infty$ and
  $A_n \norm{z}_n \le \|z\| \le 2\,D A_n \norm{z}_n$ for $ z\in
  [z_i]_n^\infty $.  Since for all $(z_i) \prec (w_i)$ we have
  $\delta_{\alpha\cdot n}([z_i], \|\cdot\| ) \le
  \ddot \delta_{\alpha\cdot n}(z_i) = \theta^n < 2 \theta ^n$,
  then  for all $(z_i) \prec (w_i)$  and all $n\in\nat$, there exists $v_n
  \in [z_i]_n^\infty$ such that $\|v_n\| \le 1$ and $\norm{v_n}_n \ge
  1/2$. Hence $A_n \le 2$, thus $\|z\| \le 4\,D \norm{z}_n$
  on $ [z_i]_n^\infty $.
 
  We have shown that for all $k\in\nat$, $\|z\| \ge (1/4\,D)\sup_{1\le n
    \le k}\norm{z}_n $ on $[w_i]_k^\infty\succ [z_i]_k^\infty$; and
  $\|z\| \le 4\,D \inf_{1\le n \le k}\norm{z}_n$ on $[z_i]_k^\infty$.
\end{proof}

We would like to directly relate the norm of an asymptotic $\ell_1$
space of bounded distortion with a norm in some Tsirelson space.
While we were unable to obtain two-sided estimates we did obtain the
following lower estimate.
\begin{prop}\label{F.h}
  Let $X$ be an asymptotic $\ell_1$ space of $D$-bounded distortion,
  $\alpha\prec \omega_1$ and suppose that $\ddot\delta_\alpha (Y) =
  \theta \in (0,1)$ for all $Y\prec X$.  Let $\ep_n\downarrow 0$.
  There exist $(w_i) \prec X$ so that for all $n$ if $w\in
  [w_i]_n^\infty$ then $\|w\| \ge (1 -\ep_n)(D+\ep_n)^{-1} \|\sum
  \|E_iw\| e_{p_i}\|_ {T(\S_\alpha,\theta-\ep_n)}$, whenever $E_1 <
  E_2 < \cdot $ are adjacent intervals, $E_i w$ denotes the
  restriction of $w$ w.r.t. $(w_i)$ and $p_i = \min E_i$.
\end{prop}

Note that the first paragraph of the proof of Theorem~\ref{F.l} shows
how to choose a subspace $X$ satisfying the above hypothesis in an
asymptotic $\ell_1$ space of bounded distortion.
\begin{proof}
  Choose $(z_i)\prec X$ so that for all $n$ there exists an equivalent
  bimonotone norm $|\cdot|_n$ on $[z_i]_n^\infty$ with $\delta_\alpha
  ((z_i)_n^\infty, |\cdot|_n) > \theta - \ep_n$.  This can be done by
  Proposition~\ref{D.na} using that $\ddot\delta_\alpha (Z) = \theta $
  for all $Z \prec X$. Hence by a diagonal argument, applying
  Corollary~\ref{E.g}, we may assume also that if $z \in \langle
  z_i\rangle_n^\infty$ and $E_1' < \cdots< E_\ell'$ are adjacent
  intervals then
  $$ |z|_n \ge (1-\ep_n) \Bigl\|\sum_1^\ell |E_i' z|_n e_{r_i}\Bigr\|_
  {T ({\S_\alpha}, \theta - \ep_n)}\ , $$ where $r_i = \min E_i' $ and
  $E_i'z$ is the restriction of $z$ w.r.t. $(z_i)$.  Using that $X$ is
  of $D$-bounded distortion and scaling $|\cdot|_n$ we may obtain
  $(w_i)\prec (z_i)$ so that for all $n$ and $w\in\langle
  w_i\rangle_n^\infty$, $\|w\|\ge |w|_n \ge \frac1{D+\ep_n} \|w\|$.
  We thus obtain for $w\in \langle w_i\rangle_n^\infty$,
\begin{eqnarray*}
  \|w\| & \ge & \left( 1-\ep_n\right) \Big\|\sum |E_i'w|_n
  e_{p_i}\Big\|_{T(\S_\alpha,\theta-\ep_n)} \\
  & \ge & {1 -\ep_n \over D+\ep_n} \Big\| \sum \|E_i'w\|
  e_{r_i}\Big\|_{T(\S_\alpha,\theta-\ep_n)}\ .
\end{eqnarray*}
 
  Now given adjacent intervals $E_1 < E_2 <\cdots$, take intervals
  $E_1' < E_2' <\cdots$ such that for all $w\in \langle
  w_i\rangle_n^\infty$, and all $i$, the restriction $E_i w$ w.r.t.
  $(w_i)$ coincides with the restriction $E_i'w$ with respect to
  $(z_i)$. Then we have $r_i = \min E_i' \ge  p_i = \min E_i$ for all $i$
  and since $\S_\alpha$ is invariant under spreading we easily get
  that $\|\sum a_i e_{r_i}\|_{T(\S_\alpha,\theta')}\ge \|\sum a_i
  e_{p_i}\|_{T(\S_\alpha,\theta')}$ for all $(a_i)$ and all $0 <
  \theta'< 1$. Thus the final lower estimate follows.
\end{proof}

The following proposition generalizes the fact that for the Tsirelson
space $T_\theta$, $D(T_\theta) \ge \theta$.
\begin{prop}\label{F.i}
  Let $X$ be an asymptotic $\ell_1$ space.  Then
$\sup \{D(Y) : Y \prec X\} \ge \sup
  \{\gamma_1^{-1} :\gamma \in\Delta (X)\}$.
\end{prop}
\begin{proof}
  Let $\gamma \in\Delta(X)$ and let $(x_i) \prec X$ $\Delta$-stabilize
  $\gamma$.  Thus for some $\ep_n\downarrow 0$, all $n$ and all
  $(y_i)\prec (x_i)_n^\infty$,
  $$\gamma_1 \ge \delta_1 (y_i) \ge \gamma_1 - \ep_n\ .$$ For
  $n\in\nat$ and $(y_i)\prec (x_i)$ define
\begin{eqnarray*}
     \lefteqn{\delta_1 (n) (y_i) = \sup \Bigl\{\delta : \|y\|\ge
  \delta\sum_1^n\|E_i y\|: y \in [y_i], Ey \mbox{ is a restriction
    w.r.t.} (y_i),}\\
  && \qquad E_1 < \cdots< E_{n_0} \mbox{ are adjacent intervals with }
      \bigcup E_i = \supp (y) \Bigr\}\ .
\end{eqnarray*}
 
  Now observe that given $\ep>0$ there exists $n_0\in\nat$  and
  $(y_i)\prec (x_i)$ so that $\delta_1(n_0)(w_i) <\gamma_1 +\ep$ for
  all $(w_i) \prec (y_i)$.  Indeed, if not, we could, by a diagonal
  argument, produce $(y_i)\prec (x_i)$ with $\delta_1(y_i) \ge
  \gamma_1+\ep$.
 
  On $[ y_i]$ define the norm
  \begin{eqnarray*}
  |y| = \sup \biggl\{ \sum_1^{n_0} \|E_iy\| &:& E_1y<\cdots < E_{n_0} y   
      \mbox{ w.r.t. } (y_i) \mbox{ and } \\
     && E_1 < \cdots < E_{n_0} \mbox{ are adjacent intervals with }
      \bigcup E_i = \supp (y)\biggr\}\ .
  \end{eqnarray*}
     
  Thus, by the choice of $(y_i)$, for all $W\prec Y = [y_i]_{\nat}$,
  there exists $w\in W$,
  $\|w\|=1$ and $|w| > \frac1{\gamma_1+\ep}$.  Also by considering
  long $\ell_1^k$-averages (see the proof of
  Proposition~\ref{B.d}) there exists $x\in W$, $\|x\|=1$ and
  $|x|<1+\ep$.  Thus $D(Y) \ge d(X,|\cdot|) \ge (1+\ep)/(\gamma_1+\ep)$.
\end{proof}

More generally, we have
\begin{prop}\label{F.j}
  Let $X$ be asymptotic $\ell_1$ and suppose that  $I_\Delta
  (X)=\alpha_0$. then
  $$\sup\{D(Y): Y\prec X\} \ge \sup \{\gamma_{\alpha_0}^{-1}
  : \gamma\in \Delta (X)\}\   .$$
\end{prop}
\begin{proof}
  We may assume $\alpha_0>1$ by Proposition~\ref{F.i}.  Thus by
  Theorem~\ref{F.c}, $\alpha_0$ is a limit ordinal.  Let
  $\alpha_n\uparrow \alpha_0$ be the ordinal sequence used in defining
  $\S_{\alpha_0}$.  Let $\gamma\in\Delta(X)$, $\ep>0$.  Then for some
  $n_0$, $\gamma_{\alpha_{n_0}} <\gamma_{\alpha_0} +\ep$.  Let $(x_i)$
  $\Delta$-stabilize $\gamma$.  Choose $(y_i)\prec (x_i)$ and  an
  equivalent norm $|\cdot|$ on $[ y_i]$  with
  $\delta_{\alpha_{n_0}}((y_i),|\cdot|)>1-\ep$.  By passing to a block
  basis of $(y_i)$ and scaling $|\cdot|$ if necessary we may assume
  that for some $D$  we have  $\|\cdot\| \le |\cdot|
  \le D\|\cdot\|$  on $[ y_i]$,
  and for all $W = [w_i] \prec Y$ there exists $w\in W$,
  $\|w\|=1$ and $|w|<1+\ep$.  Since $\gamma_{\alpha_{n_0}}
  <\gamma_{\alpha_0} +\ep$ there exists $z\in W$ with $\|z\|=1$ and
  $\sum_1^\ell \|z_i \| \ge 1/(\gamma_{\alpha_0} +\ep)$,
  for some decomposition $z = \sum_1^\ell z_i $  where $(z_i)_1^\ell$
  is $\alpha_{n_0}$-admissible w.r.t. $(w_i)$.
  Hence $|z| \ge (1 - \ep) \sum |z_i| \ge (1 - \ep) \sum \|z_i\| \ge
  (1-\ep)/(\gamma_{\alpha_0}+\ep)$. Comparing the norms $|z|$ and
  $\|z\|$ we get  $D(Y,|\cdot|) >  (1 -\ep)(1+\ep)/(\gamma_
  {\alpha_0}+\ep)$.
\end{proof}

We have a simple corollary.

\begin{cor}\label{F.k}
  Let $X$ be asymptotic $\ell_1$ with $I_\Delta (X) = I_\Delta (Y) =
  \alpha_0$ for all $Y\prec X$.  If $\ddot\delta_{\alpha_0} (X)=0$
  then no subspace of $X$ is of bounded distortion.
\end{cor}

\medskip\noindent
Edward Odell:
Department of Mathematics,
The University of Texas at Austin,\\
Austin, TX 78712, USA \\
e-mail: {\small\tt%
  odell@math.utexas.edu}

\medskip\noindent
Nicole Tomczak-Jaegermann:
Department of Mathematical Sciences,
University of Alberta,\\
Edmonton, Alberta, T6G 2G1, Canada\\
e-mail: {\small\tt%
  ntomczak@vega.math.ualberta.ca}

\medskip\noindent
Roy Wagner:
School of Mathematical Sciences, Sackler Faculty of Exact Sciences,\\
Tel Aviv University, Tel Aviv 69978, Israel\\
e-mail: {\small\tt%
      pasolini@math.tau.ac.il}


\begin{thebibliography}{99}
\frenchspacing

\bibitem{AA} D. Alspach and S. Argyros:  {\em Complexity of weakly
    null sequences}, Diss. Math., {\bf 321} (1992), 1--44. 

\bibitem{AO} G. Androulakis and E. Odell:  {\em Distorting mixed
    Tsirelson spaces}, preprint.
 
\bibitem{AD} S. Argyros and I. Deliyanni:  {\em Examples of
    asymptotically $\ell^1$ Banach spaces}, preprint.

\bibitem{AMT} S. Argyros, S. Merkourakis and A. Tsarpolias:  {\em
    Convex unconditionality and sumability of  weakly null sequences},
    preprint.

\bibitem{Belle} S. Bellenot: {\em The Banach space $T$ and the fast
    growing  hierarchy from logic}, Israel J.~of Math., {\bf 47}
  (1984), 305--313.
 
\bibitem{B} J. Bourgain:  {\em On convergent sequences of continuous
    functions}, Bull. Soc. Math. Bel., {\bf 32} (1980), 235--249.
 
\bibitem{CJT} P. G. Casazza, W. B. Johnson and L. Tzafriri:  {\em On
    Tsirelson's space}, Israel J. Math., {\bf 47} (1984), 81--98.

\bibitem{CS} P. G. Casazza and T. J. Shura:  {\em Tsirelson's Space},
  Lecture Notes in Math., vol.~1363, Springer-Verlag, Berlin and New
  York, 1989.
 
\bibitem{DK} I. Deliyanni and D. Kutzarova:  {\em On some asymptotic
    $\ell_1$ Banach spaces}, preprint.

\bibitem{FJ} T. Figiel and W. B. Johnson:  {\em A uniformly convex
    Banach space which contains no $\ell_p$}, Comp. Math., {\bf29}
  (1974), 179--190.

\bibitem{Jam} R. C. James: {\em  Uniformly non-square Banach spaces},
     Ann. of Math., {\bf80} (1964), 542--550.

\bibitem{Judd} R. Judd: in preparation. 
 
\bibitem{JO} R. Judd and E. Odell:  {\em Concerning the Bourgain
    $\ell_1$ index of a Banach space}, preprint.

\bibitem{Ku} K. Kuratowski: {\em Topology}, vol I, academic Press,
  New York 1966.

\bibitem{Leung} D. Leung, {\em Compact subsets of ${\cal P}_{<\infty} (\nat)$ 
with applications to the embedding of symmetric sequence spaces into 
$C(\alpha)$}, preprint.

\bibitem{LT} J. Lindenstrauss and L. Tzafriri: {\em Classical Banach
    spaces}, vol I, Springer Verlag 1977.

\bibitem{Ma} B. Maurey:  {\em A remark about distortion}, Operator
  Theory: Advances and Applications, {\bf77} (1995), 131--142.
 
\bibitem{MMT} B. Maurey, V.D. Milman and N. Tomczak-Jaegermann:  {\em
    Asymptotic infinite-dimensional theory of Banach spaces}, Operator
  Theory: Advances and Applications, {\bf77} (1995), 149--175.

\bibitem{MiT} V. Milman and N. Tomczak-Jaegermann:  {\em Asymptotic
    $\ell_p$ spaces and bounded distortion}, (Bor-Luh Lin and W.B.
  Johnson, eds.), Contemporary Math., {\bf144} Amer. Math. Soc.,(1993),
  173--195.
 
\bibitem{M} J. D. Monk:  {\em Introduction to Set Theory}, McGraw-Hill,
  1969.
 
\bibitem{OS1} E. Odell and Th. Schlumprecht:  {\em The distortion
    problem}, GAFA, {\bf 3} (1993), 201--207.

\bibitem{OS2} E. Odell and Th. Schlumprecht:  {\em The distortion
    problem}, Acta Math., {\bf173} (1994), 259--281.
 
\bibitem{R2} H. Rosenthal:  {\em A characterization of Banach spaces
    containing $\ell_1$}, Proc. Nat. Acad. Sci. (U.S.A.) {\bf71}
  (1974), 2411--2413.
 
\bibitem{To} N. Tomczak-Jaegermann:  {\em Banach spaces of type $p >1$
    have arbitrarily distortable subspaces}, preprint. 

\bibitem{Tom} N. Tomczak-Jaegermann: Notes privately distributed,
  December 1994.

\bibitem{Ts} B. S. Tsirelson:  {\em Not every Banach space contains
    $\ell_p$ or $c_0$}, Funct. Anal. Appl., {\bf8} (1974), 138--141.


\end{thebibliography}
\end{document}